\documentclass{IEEEtran}

 \pdfoutput=1

\usepackage{graphicx}
\usepackage{amsmath}
\usepackage{amssymb}
\usepackage{amsfonts}
\usepackage{mathrsfs}
\usepackage{dsfont}
\makeatletter
\def\amsbb{\use@mathgroup \M@U \symAMSb}
\makeatother
\usepackage{color}
\usepackage{subfigure}
\usepackage[usenames,dvipsnames]{xcolor}

\newtheorem{theorem}{Theorem}
\newtheorem{remark}{Remark}
\newtheorem{definition}{Definition}

\newtheorem{proposition}{Proposition}
\newtheorem{corollary}{Corollary}
\newtheorem{assumption}{Assumption}

\usepackage{balance}

\newcommand{\norm}[1]{\left\lVert#1\right\rVert}
\usepackage{mathrsfs}

\title{\color{black}Primary frequency regulation with load-side participation: stability and optimality}

\author{Andreas Kasis\thanks{Andreas Kasis and Ioannis Lestas are with the Department of Engineering, University of Cambridge, Trumpington Street, Cambridge, CB2 1PZ, United Kingdom; e-mails: ak647@cam.ac.uk, icl20@cam.ac.uk}, Eoin Devane\thanks{Eoin Devane is with the Cambridge Centre for Analysis, Centre for Mathematical Sciences, University of Cambridge, Wilberforce Road, Cambridge, CB3 0WA, United Kingdom; e-mail: esmd2@cam.ac.uk}, and Ioannis Lestas}

\begin{document}

\maketitle

\begin{abstract}
We present a method to design distributed generation and demand control schemes for primary frequency regulation in power networks that guarantee asymptotic stability and ensure fairness of allocation. We impose a passivity condition on net power supply variables and provide explicit steady state conditions on a general class of generation and demand control dynamics that ensure convergence of solutions to equilibria that solve an appropriately constructed network optimization problem. We also show that the inclusion of controllable demand results in a drop in steady state frequency deviations. We discuss how various classes of dynamics used in recent studies fit within our framework and show that this allows for less conservative stability and optimality conditions. We illustrate our results with simulations on the IEEE 68 bus system and observe that both static and dynamic demand response schemes that fit within our framework offer improved transient and steady state behavior compared with control of generation alone. The dynamic scheme is also seen to enhance the robustness of the system to time-delays.
\end{abstract}

\section{Introduction}

Large scale integration of renewable sources of energy within the power grid is expected to cause fast changes in generation, making power imbalances increasingly frequent due to the inability of conventional means of generation to counter-balance them. Load participation is considered to be one potential solution to this problem, providing fast response to power changes. Household appliances like air conditioning units, heaters, and refrigerators can be controlled to adjust frequency and regulate power imbalances. Although the idea dates back to the 1970s~\cite{Schweppe}, research attention has recently increasingly focused on the concept of controllable demand~\cite{Recent Load Control}, with particular consideration given to its use for primary control~\cite{Trudnowski},~\cite{Primary2}. Primary or droop control operates on a timescale of tens of seconds and is responsible for ensuring that generation and demand are always balanced~\cite{Bergen_Vittal}. It functions by adjusting governor set-points based on local frequency deviations and is completely decentralized. Simulation studies and field tests (e.g.~\cite{Pacific NW LAB,LIPA}) have shown the feasibility of including controllable demand in the droop control process. Furthermore, they have shown that local frequency measurements are sufficient as a control signal for load participation in primary control.
%

An issue of fairness is raised if appliances are to be used as a means of controllable demand. Recent studies have attempted to address this issue by devising control schemes which solve an optimization problem guaranteeing fair allocation between them. This approach has been studied for primary and also for secondary control. We consider here primary rather than secondary control in order to avoid the additional communication that would be necessary to get a fair allocation if controllable demand were used in the latter case. This is because it is evident that a synchronizing variable is necessary to achieve optimality, allowing all nodes to adapt their generation and controllable demand so as to attain equal marginal costs. In primary control, frequency deviation from the nominal value can be used for this purpose, allowing decentralized control to be achieved~\cite{P19,P34 - Low}. {\color{black}In ~\cite{P19} only primary control in loads and linearised power flows are considered while in~\cite{P34 - Low} these assumptions are waived, considering also generation control and non-linear power flows.} In secondary control, where frequency deviations return to zero, a different variable needs to be used for synchronization. A common approach is to introduce additional information exchange between neighbouring buses in order to achieve this~\cite{LC2,LC3,LC4,trip} {\color{black}or have special network structures like star topology, as in~\cite{P52}.}


{\color{black}In this paper}, we consider a network model described by nonlinear swing equations. We consider a general class of dynamics for power generation and controllable demand, on which we impose appropriate conditions so as to achieve stability of the equilibrium points and an optimization interpretation of those. This allows us to guarantee, for a wide variety of possible generation and demand dynamics, convergence to a power allocation that solves an appropriately constructed optimization problem, thus ensuring fairness in this allocation. The class of dynamics considered incorporates control schemes using only local frequency measurements as input signals, and we demonstrate that this is sufficient to enable them to take the right decisions so as to converge to a global optimum, thus allowing distributed control. Further, we also prove that the inclusion of controllable demand can aid in frequency control by reducing the steady state frequency error. We illustrate the applicability of our approach by demonstrating that various dynamics that have been used in recent interesting studies, such as~\cite{P19} and~\cite{P34 - Low}, can be incorporated within our framework, and we show that the analysis presented in the paper allows for less conservative stability and optimality conditions.

{\color{black}
It should be noted that one of the distinctive features of our analysis is that optimality of the power allocation is provided via appropriate conditions on the input/output properties of the systems considered.
In our companion paper \cite{paper2} we show how additional local information can be exploited to relax a passivity condition used in this paper to deduce convergence.}

The paper is organized as follows. {\color{black} In Sections~\ref{Notation} and~\ref{sec: Preliminaries} we {\color{black}start with some basic notation and definitions} that are used in the rest of the paper.} In Section~\ref{Network_Model}, we present the power network model and introduce the class of generation and load dynamics that will be considered. Section~\ref{Network_Model} introduces the optimization problem that we consider and Section~\ref{Main_Results} presents our main results, which are proved in Appendix A. In Section~\ref{Discussion}, we discuss how our analysis relates to other important studies. Section~\ref{Simulation} illustrates our results through simulations on the IEEE 68 bus system. Finally, conclusions are drawn in Section~\ref{Discussion_Future_Work}. {\color{black}In Appendix B we present a generalization of the optimality results in the paper whereby the cost functions considered are allowed to have discontinuous derivatives, which, as discussed in the paper, can be important in practical implementations.}

\section{Notation}\label{Notation}

{\color{black}
{\color{black} Real numbers are denoted by $\mathbb{R}$, and the set of n-dimensional vectors with real entries is denoted by $\mathbb{R}^n$.}
{\color{black}For a function $f(q)$, we denote its first derivative by $f'(q) = \tfrac{d}{dq} f(q)$.} The expression $f^{-1}(w)$ represents the preimage of the point $w$ under the function $f$, i.e. $f^{-1}(w) = \{q \colon f(q) = w\}$. When the function $f$ is invertible, $f^{-1}$ then defines the inverse function of $f$. {\color{black}{\color{black}
A function $f:\mathbb{R}^n \rightarrow \mathbb{R}$
is said to be positive definite on a neighbourhood D around the origin}
if $f(0) = 0$  and $f(x) > 0$ for every non-zero $x \in D$. {\color{black}It is positive semidefinite if the {\color{black}inequality $>0$ is replaced by
$\ge0$}.}} {\color{black}Furthermore, for $a,b \in \mathbb{R}$, $a\leq b$, the expression $[q]^b_a$ will be used to denote $\max\{\min\{q,b\},a\}$.} The indicator function $\mathds{1}_S : \mathbb{R}^n \rightarrow \{0,1\}$ of a set $S \subseteq \mathbb{R}^n$ takes the {\color{black}value~$1$ if its argument belongs to the set $S$ and $0$} otherwise.

{\color{black} For a function of time {\color{black}$q(t)$}, {\color{black}we denote its derivative with respect to time by $\dot{q}$ 
and its Laplace transform by} $\hat{q} = \int_0^\infty e^{-st} q(t) \, dt$.
{\color{black}We use $\mathcal{L}_2$} to denote} the Hilbert space of functions {\color{black}$f:\mathbb{R} \rightarrow \mathbb{R}$} with finite norm
{\color{black}$\norm{f}_2$ = $\sqrt{\int_{0}^{\infty} f^2(t)dt}$.}
{\color{black}For a system as in~\eqref{dynsys} where $x=\bar x$, $u=y=0$ is an equilibrium point, its} $\mathcal{L}_2$-gain is defined as {\color{black} {\color{black}the supremum of the ratio of the $\mathcal{L}_2$-norms of the output $y$ and the input $u$, i.e. $\sup_{\norm{u}_2\neq0}\frac{\norm{y}_2}{\norm{u}_2}$, with $x(0)=\bar x$. For a stable linear system with transfer function $T(s)$ it is known that its $\mathcal{L}_2$-gain is given by its $\infty$-norm $\norm{T}_{\infty} := \sup\limits_{\omega} |T(j\omega)|$}.}
}

\section{Preliminaries} \label{sec: Preliminaries}

{\color{black} Throughout the paper we will consider dynamical systems
with input $u(t) \in \mathbb{R}$, state $x(t) \in \mathbb{R}^m$,}
and output $y(t) \in \mathbb{R}$ with a state space realization of the form
\begin{equation} \label{dynsys}
\begin{aligned}
&\dot{x} = f(x,u),\\
&y = g(x,u),
\end{aligned}
\end{equation}
where $f : \mathbb{R}^m \times \mathbb{R} \to \mathbb{R}^m$ is locally Lipschitz and $g : \mathbb{R}^m \times \mathbb{R} \to \mathbb{R}$ is continuous. {\color{black}We assume in system~\eqref{dynsys} 
that given} any constant input $u(t) \equiv \bar{u} \in \mathbb{R}$, there exists a {\color{black}unique locally asymptotically stable equilibrium point $\bar{x} \in \mathbb{R}^m$, {\color{black}i.e.} $f(\bar{x}, \bar{u}) = 0$. {\color{black} The region of attraction\footnote{\color{black}{\color{black}That is}, for the constant input $u = \bar{u}$, any solution $x(t)$ of~\eqref{dynsys} with initial condition $x(0) \in X_0$ must satisfy $x(t) \to \bar{x}$ as $t \to \infty$. The definition of local asymptotic stability also implies that this region $X_0$ has nonempty interior.} {\color{black}of~$\bar{x}$} is denoted by $X_0$. {\color{black}We} also define the static input-state characteristic {\color{black} map $k_x : \mathbb{R} \to \mathbb{R}^m$,}}
\begin{equation} \label{ischar}
{\color{black}k_x(\bar{u}) := \bar{x}.}
\end{equation}
Based on this, we can also define the static input-output characteristic map {\color{black}
{
$k_y : \mathbb{R} \to \mathbb{R}$, }}}
\begin{equation}
{\color{black}k_y(\bar{u}) := g(k_x(\bar{u}), \bar{u}).}
\label{iochar}
\end{equation}
{\color{black}The requirement that {\color{black}for any constant input to~\eqref{dynsys} there exists a unique equilibrium point}, could be relaxed to require only isolated {\color{black}equilibrium points}, however, we assume it here to simplify the presentation.
}

\section{Problem formulation} \label{Network_Model}

\subsection{Network model}

The power network model is described by a graph $(N,E)$ where $N = \{1,2,\dots,|N|\}$ is the set of buses and $E \subseteq N \times N$ the set of transmission lines connecting the buses. There are two types of buses in the network, generation and load buses. Their main difference is that generation buses have non-zero generation inertia in contrast with load buses. Correspondingly, only generation buses have nontrivial generation dynamics. Let $G = \{1,2,\dots,|G|\}$ and $L= \{|G|+1,\dots,|N|\}$ be the sets of generation and load buses such that $|G| + |L| = |N|$. Furthermore, we use $(i,j)$ to denote the link connecting buses $i$ and $j$ and assume that the graph $(N,E)$ is directed with arbitrary direction, so that if $(i,j) \in E$ then $(j,i) \notin E$. For each $j \in N$, we use $i:i\rightarrow j$ and $k:j\rightarrow k$ to denote the sets of buses that are predecessors and successors of bus~$j$ respectively. It is important to note that the form of the dynamics in~\eqref{sys1}--\eqref{sys2} below is unaltered by any change in the graph ordering, and all of our results are independent of the choice of direction. We also assume that $(N,E)$ is connected. The following assumptions are made for the network: \newline
1) Bus voltage magnitudes are $|V_j| = 1$ p.u. for all $j \in N$. \newline
2) Lines $(i,j) \in E$ are lossless and characterized by their susceptances $B_{ij} = B_{ji} > 0$. \newline
3) Reactive power flows do not affect bus voltage phase angles and frequencies.

The rate of change of frequency at the generation buses can then be represented using swing equations, while power must be conserved at each of the load buses. This motivates the following system dynamics (e.g.~\cite{Bergen_Vittal})
\begin{subequations} \label{sys1}
\begin{equation}
\dot{\eta}_{ij} = \omega_i - \omega_j, \; (i,j) \in E, \label{sys1a}
\end{equation}
\begin{equation}
M_j \dot{\omega}_j = - p_j ^L + p_j^M - (d^c_j + d^u_j) - \sum_{k:j\rightarrow k} p_{jk} + \sum_{i:i\rightarrow j} p_{ij}, \; j\in G, \label{sys1b}
\end{equation}
\begin{equation}
 0 = - p_j ^L - (d^c_j + d^u_j) - \sum_{k:j\rightarrow k} p_{jk} + \sum_{i:i\rightarrow j} p_{ij}, \; j\in L, \label{sys1c}
\end{equation}
\begin{equation}
{\color{black} p_{ij}=B_{ij} \sin\eta _{ij} - p_{ij}^{nom}, \; (i,j) \in E.} \label{sys1d}
\end{equation}
\end{subequations}
In system~\eqref{sys1} the quantities  $p^M_j$, $\omega_j$, $d^c_j$ and $d^u_j$ are time-dependent variables representing, respectively, deviations from a nominal value\footnote{{\color{black}A nominal value of a variable is defined as its value at an equilibrium of \eqref{sys1} with frequency equal to the nominal value of 50Hz (or 60Hz).}} for the mechanical power injection to the generator bus $j$, and the frequency, controllable load and uncontrollable frequency-dependent load present at any bus~$j$.
{\color{black}{\color{black}The variables}
 $\eta_{ij}$ and $p_{ij}$ {\color{black}represent, respectively,} the power angle difference\footnote{The quantities $\eta_{ij}$ represent the phase differences between buses $i$ and $j$, given by $\theta_i - \theta_j$. The angles themselves must also satisfy $\dot{\theta}_j = \omega_j$ at all $j \in N$, however, we omit this equation in \eqref{sys1} since the power transfers $p$ are functions only of the phase differences.} and  the deviation from nominal value $p_{ij}^{nom}$ for the} power transmitted from bus $i$ to bus $j$.
The constant $M_j > 0$ denotes the generator inertia.
{\color{black}We shall study the response of system~\eqref{sys1} to a step change $p_j^L$ in the uncontrollable demand at each bus $j$.} To investigate decentralized control schemes for {\color{black}generation and controllable load} based upon local measurements of the frequency alone, we close the loop in~\eqref{sys1} by determining each of $p^M_j$, $d^c_j$, and $d^u_j$ as outputs from independent systems of the {\color{black}form in} Section~\ref{sec: Preliminaries} with inputs given by the negative of the local frequency,
\begin{subequations} \label{sys2}
\begin{equation} \label{sys2p}
\begin{aligned}
&\dot{x}^{M,j} = f^{M,j}(x^{M,j},-\omega_j), \\
&p^M_j = g^{M,j}(x^{M,j},-\omega_j),
\end{aligned} \hspace{2em}j \in G,
\end{equation}
\begin{equation} \label{sys2dc}
\begin{aligned}
&\dot{x}^{c,j} = f^{c,j}(x^{c,j},-\omega_j), \\
&d^c_j = g^{c,j}(x^{c,j},-\omega_j),
\end{aligned} \hspace{2em}j \in N,
\end{equation}
\begin{equation} \label{sys2du}
\begin{aligned}
&\dot{x}^{u,j} = f^{u,j}(x^{u,j},-\omega_j), \\
&d^u_j = g^{u,j}(x^{u,j},-\omega_j),
\end{aligned} \hspace{2em}j \in N.
\end{equation}
\end{subequations}
{\color{black}For convenience in the notation}, we collect\footnote{{\color{black}Note that} each local variable (e.g. $x^{M,j}$) is a vector with multiple components.} the variables in~\eqref{sys2} into the vectors $x^M = [x^{M,j}]_{j \in G}$, $x^c = [x^{c,j}]_{j \in N}$, and $x^u = [x^{u,j}]_{j \in N}$. These quantities represent the internal states of the dynamical systems used to update the desired outputs $p^M_j$, $d^c_j$, and $d^u_j$. The variables $p^M_j$ and $d^c_j$ are controllable, so we have freedom in our analysis to design certain properties of the dynamics in~\eqref{sys2p} and~\eqref{sys2dc}. By contrast, $d^u_j$ represents uncontrollable load and the dynamics
in~\eqref{sys2du} are thus fixed. {\color{black}It should be noted that the systems in~\eqref{sys2} can be heterogeneous and of arbitrary dimension.

Throughout the paper we aim to characterize broad classes of 
 dynamics associated with generation and demand, so that
  {\em stability} and {\em optimality} can be guaranteed for the equilibrium points of the overall interconnected system~\eqref{sys1}--\eqref{sys2}. Various examples of those will be discussed in Section~\ref{Discussion}, with simulations also provided in section \ref{Simulation}.
}



\subsection{Equilibrium analysis}

We now quantify what is meant by an equilibrium of the interconnected system~\eqref{sys1}--\eqref{sys2}.

\begin{definition} \label{eqbrdef}
The constants\footnote{By constant {\color{black}we mean a variable independent} of time.} $(\eta^*, \omega^*, x^{M,*}, x^{c,*}, x^{u,*})$ define an equilibrium of the system~\eqref{sys1}--\eqref{sys2} if the following hold
\begin{subequations} \label{eqbr}
\begin{equation}
0 = \omega^*_i - \omega^*_j, \; (i,j) \in E, \label{eqbr1}
\end{equation}
\begin{equation}
0 = - p_j^L + p_j^{M,*} - (d^{c,*}_j + d^{u,*}_j) - \sum_{k:j\rightarrow k} p^*_{jk} + \sum_{i:i\rightarrow j} p^*_{ij}, \; j\in G, \label{eqbr2}
\end{equation}
\begin{equation}
0 = - p_j^L - (d^{c,*}_j + d^{u,*}_j) - \sum_{k:j\rightarrow k} p^*_{jk} + \sum_{i:i\rightarrow j} p^*_{ij}, \; j\in L, \label{eqbr3}
\end{equation}
\begin{equation}
x^{M,j,*} = k_{x^{M,j}} (-\omega^*_j), \; j \in G, \label{eqbr8}
\end{equation}
\begin{equation}
x^{c,j,*} = k_{x^{c,j}} (-\omega^*_j), \; j \in N, \label{eqbr9}
\end{equation}
\begin{equation}
x^{u,j,*} = k_{x^{u,j}} (-\omega^*_j), \; j \in N, \label{eqbr10}
\end{equation}
where the quantities in~\eqref{eqbr2} and~\eqref{eqbr3} are given by
\begin{equation}
{\color{black} p^*_{ij}=B_{ij} \sin \eta^*_{ij} - p_{ij}^{nom}, \; (i,j) \in E, }\label{eqbr4}
\end{equation}
\begin{equation}
p^{M,*}_j = k_{p^M_j} (-\omega_j^*), \; j \in G, \label{eqbr5}
\end{equation}
\begin{equation}
d^{c,*}_j = k_{d^c_j} (-\omega_j^*), \; j \in N, \label{eqbr6}
\end{equation}
\begin{equation}
d^{u,*}_j = k_{d^u_j} (-\omega_j^*), \; j \in N. \label{eqbr7}
\end{equation}
\end{subequations}
We call~\eqref{eqbr} the equilibrium conditions for the system~\eqref{sys1}--\eqref{sys2}.
\end{definition}

\begin{remark}
For any equilibrium with a given frequency value~$\omega^*$, the uniqueness in the definition of the static input-state characteristic in Section~\ref{sec: Preliminaries} immediately shows that the values of $p^{M,*}$, $d^{c,*}$, and $d^{c,*}$ are all guaranteed to be unique. By contrast, there can in general be multiple choices of $\eta^*$ and~$p^*$ such that the equilibrium equations~\eqref{eqbr} remain valid, and thus it follows that the equilibrium power transfers $p^*_{ij}$ between regions need not be unique. It can be shown that they become unique under prescribed conditions on the network structure, such as when $(N,E)$ has tree topology.
\end{remark}

Note that the static input-output characteristic maps $k_{p^M_j}$, $k_{d^c_j}$, and $k_{d^u_j}$, relating power generation/demand with frequency, as defined in~\eqref{iochar}, completely characterize the effect of the dynamics~\eqref{sys2} on the behavior of the power system~\eqref{sys1} at equilibrium. In our analysis, we will consider a class of dynamics within~\eqref{sys2} for which the set of equilibrium points in Definition~\ref{eqbrdef} that satisfy Assumption~\ref{assum1} below is asymptotically stable. Within this class, we then consider appropriate conditions on these characteristic maps such that the values of the variables defined in~\eqref{eqbr5}--\eqref{eqbr7} are {\em optimal} for an appropriately constructed network optimization problem.

Throughout the remainder of the paper we suppose that there exists some equilibrium of~\eqref{sys1}--\eqref{sys2} as defined in Definition~\ref{eqbrdef}. We let $(\eta^*, \omega^*, x^{M,*}, x^{c,*}, x^{u,*})$ denote any such equilibrium and use $(p^*, p^{M,*}, d^{c,*}, d^{u,*})$ to represent the corresponding quantities defined in~\eqref{eqbr4}--\eqref{eqbr7}. We now impose an assumption on the equilibrium, which can be~interpreted as a security constraint for the power flows generated.
{\color{black}
\begin{assumption} \label{assum1}
$| \eta^*_{ij} | < \tfrac{\pi}{2}$ for all $(i,j) \in E$.
\end{assumption}
The stability and optimality properties of such an equilibrium point will now be studied in the sections that follow.
}


\subsection{Combined passive dynamics from generation and load dynamics } \label{Passive Dynamics}

In terms of the outputs from~\eqref{sys2}, define the net supply variables
\begin{subequations} \label{ssys}
\begin{equation} \label{ssysG}
{\color{black}s^G_j = p^M_j - (d^c_j + d^u_j), \; j \in G,}
\end{equation}
\begin{equation} \label{ssysL}
{\color{black}s^{L}_j = -(d^c_j + d^u_j), \; j \in L.}
\end{equation}
\end{subequations}
Correspondingly, their values at equilibrium can be written as $s^{G,*}_j = p^{M,*}_j - (d^{c,*}_j + d^{u,*}_j)$ and $s^{L,*}_j = - (d^{c,*}_j + d^{u,*}_j)$.

The variables defined {\color{black}in~\eqref{ssys}} evolve according to the dynamics described in~\eqref{sys2}. Consequently, $s^G_j$ and $s^L_j$ can be viewed as outputs from these combined dynamical systems with inputs~$-\omega_j$.

{\color{black}We now introduce a notion of passivity for systems of the form~\eqref{dynsys} {\color{black}which we will use for the dynamics of the supply variables defined in~\eqref{ssys} to prove our main stability results.}}

{\color{black}
\begin{definition}\label{Passivity_Definition}
The system~\eqref{dynsys} is said to be locally input strictly passive {\color{black}about the constant input values $\bar{u}$ and the constant state values~$\bar{x}$} if {\color{black}there exist open neighbourhoods $U$ of $\bar{u}$ and $X$ of $\bar{x}$ and a continuously differentiable, positive semidefinite function $V(x)$} (called the storage function) such {\color{black}that, for all $u \in U$ and all $x \in X$,}
\begin{equation}\label{In_passivity_def}
\dot{V}(x) \leq (u - \bar{u})^T (y - \bar{y}) - \phi(u - \bar{u})
\end{equation}
where $\phi$ is a positive definite function and $\bar{y} = k_y (\bar{u})$.
\end{definition}
}

We now suppose that {\color{black}the} maps from negative frequency to supply satisfy {\color{black}this condition of input strict} passivity about {\color{black}equilibrium. Note that this is a decentralized condition.

{\color{black}
\begin{assumption} \label{assum2}
{\color{black}Each of the} systems defined in~\eqref{sys2} with inputs $-\omega_j$ and outputs given by~{\color{black}\eqref{ssysG} and~\eqref{ssysL} respectively} are locally input strictly passive {\color{black}about their equilibrium values $-\omega^*_j$ and {\color{black}$(x^{M,j,*}, x^{c,j,*}, x^{u,j,*})$}, in the sense {\color{black} described} in} Definition~\ref{Passivity_Definition}{\color{black}.}
\end{assumption}
}}

%

\begin{remark}
Note that, in Assumption~\ref{assum2}, we assume only a passivity property for the systems~\eqref{sys2} without specifying the precise form of the systems, which permits the inclusion of a broad class of generation and load dynamics. Moreover, the fact {\color{black}that passivity} is assumed only for the net supply dynamics, rather than for the generation and load dynamics individually, can permit the analysis of systems incorporating dynamics that are not individually passive. We will see in {\color{black}Section~\ref{Discussion}} that various classes of generation and load dynamics that have been investigated in the literature fit within the present framework.
\end{remark}

{\color{black} It should be noted that for linear systems input strict passivity can easily be checked from strict positive realness of the transfer function or numerically using the KYP Lemma \cite{Khalil}. It holds for sign preserving  static nonlinearities, first-order linear dynamics, and also higher-order generation/load dynamics with sufficiently large damping, as it will be discussed in Section~\ref{Discussion}.}

\subsection{Optimal supply and load control} \label{sec:optim}

We aim to explore how the generated power and controllable loads may be adjusted to meet the step change $p^L$ in frequency-independent load in a way that minimizes the total cost that comes from the extra power generated and the disutility of loads. We now introduce an optimization problem, which we call the optimal supply and load control problem (OSLC), that can be used to achieve this goal.

{\color{black}
Suppose that costs $C_j(p^M_j)$ and $C_{dj}(d^c_j)$ are incurred for deviations $p^M_j$
and $d^c_j$ in generation and controllable load respectively.}
Furthermore, some additional cost is incurred due to any change in frequency which alters the demand from uncontrollable frequency-dependent loads. We represent this by an integral cost in terms of a function $h_j$ which is explicitly determined by the dynamics in~\eqref{sys2du} as
\begin{equation}
h_j(z) = k_{d^u_j}(-z) \text{ for all } \bar{z} \in \mathbb{R}. \label{hdef}
\end{equation}
The total cost within OSLC then sums all the above costs, and the problem is to choose the vectors $p^M$, $d^c$, and $d^u$  that minimize this total cost and simultaneously achieve power balance, while satisfying physical {\color{black}saturation {\color{black}constraints}.
\begin{equation}
\begin{aligned}
&\hspace{2em}\underline{\text{OSLC:}} \\
&\min_{p^M,d^c,d^u} \sum\limits_{j\in G} C_{j} (p_{j}^M) + \sum\limits_{j\in N} \Big( C_{dj} (d^c_{j}) + \int_0^{d^u_j} \! h_{j}^{-1}(z) \, dz \Big) \hspace{-1.5em}\\
&\text{subject to } \sum\limits_{j\in  G} p_j^M = \sum\limits_{j\in  N} (d^c_j + d^u_j +p_j^L), \\
& p^{M,min}_j \leq p^M_j \leq p^{M,max}_j \,,\, \forall j \in G, \\
& d^{c,min}_j \leq d^c_j \leq d^{c,max}_j \,,\, \forall j \in N,
 \label{Problem_To_Min}
\end{aligned}
\end{equation}
{\color{black}where} $p^{M,min}_j,p^{M,max}_j, d^{c,min}_j$, and $d^{c,max}_j$ are the bounds for generation and controllable demand {\color{black}respectively at bus~$j$}. The equality constraint {\color{black}in~\eqref{Problem_To_Min}} represents conservation of power by specifying that all the extra frequency-independent load is matched by the total additional generation plus all the deviations in frequency-dependent loads.}

\begin{remark}
The variables $p^M$ and $d^c$ within~\eqref{Problem_To_Min} represent the variables that can be directly controlled (generation and controllable demand), while the variable $d^u$ can be controlled only indirectly by effecting changes in the system frequencies (uncontrollable frequency-dependent demand). Therefore, we aim to specify properties of the control dynamics in~\eqref{sys2p}--\eqref{sys2dc} that ensure that the quantities $p^M$ and $d^c$, along with the system frequencies, converge to values at which optimality in~\eqref{Problem_To_Min} can be {\color{black}guaranteed.}
\end{remark}

{\color{black}\subsection{Additional conditions}}

To guarantee convergence and optimality, we will {\color{black}require additional conditions} on the behavior of the systems~\eqref{sys1}--\eqref{sys2} and the structure of the optimization problem~\eqref{Problem_To_Min}. {\color{black}The assumptions introduced are all of practical relevance, and we will see in Section~\ref{Discussion} that the framework considered encompasses a number of important examples that have been investigated within the literature.}

The first two of these conditions are needed for our proof of the convergence result in Theorem~\ref{convthm}. Within the second of these we will denote $\omega^G = [\omega_j]_{j \in G}$ and $\omega^L = [\omega_j]_{j \in L}$.

\begin{assumption} \label{assum3}
{\color{black}The storage functions in Assumption~\ref{assum2}} have strict local minima at the points $(x^{M,j,*},x^{c,j,*},x^{u,j,*})$ and $(x^{c,j,*},x^{u,j,*})$ respectively.
\end{assumption}

{\color{black}
\begin{remark}
In practice, Assumption~\ref{assum3} is often trivially satisfied. For instance, whenever the vector fields in~\eqref{sys2} are continuously differentiable, then by {\color{black} considering a linearization about equilibrium,} 
the KYP Lemma can be seen to generate a storage function satisfying Assumption~\ref{assum3}, {\color{black} when the linearized system is controllable and observable}.
\end{remark}
}

{\color{black}
\begin{assumption} \label{assum4}
There exists an open neighbourhood~$T$ of $(\eta^*, \omega^{G,*}, x^{M,*},x^{c,*},x^{u,*})$ such that at any time instant $t$, $\omega^L(t)$ is uniquely determined by the system states $(\eta(t), \omega^G(t),$ $x^M(t), x^c(t), x^u(t)) \in T$ and equations \eqref{sys1}--\eqref{sys2}.
\end{assumption}
}

\begin{remark}
Assumption~\ref{assum4} is a technical assumption that is required in order for the system~\eqref{sys1}--\eqref{sys2} to have a locally well-defined state space realization. This is needed in order to apply Lasalle's Theorem to analyze stability in the proof of Theorem~\ref{convthm} {\color{black}below}. Without Assumption~\ref{assum4}, stability could be analyzed through more technical approaches such as the singular perturbation analysis discussed in~\cite[Section 6.4]{sastry}.
\end{remark}
\begin{remark}
Assumption~\ref{assum4} can often be verified by using the Implicit Function Theorem to generate decentralized algebraic conditions under which it is guaranteed to hold. {\color{black}For instance, Assumption~\ref{assum4} always holds if
in~\eqref{sys2} we have for all $j \in L$
\begin{equation} \label{a4cond1}
\frac{\partial g^{c,j}}{\partial \omega_j}
+ \frac{\partial g^{u,j}}{\partial \omega_j} \neq 0,
\end{equation}
at the equilibrium point. If the functions $g^{c,j}$ and $g^{u,j}$ have no explicit dependence on $\omega_j$ 
 satisfying  the following condition at the equilibrium point is also sufficient
\begin{align} \label{a4cond2}
&\sum_i \frac{\partial g^{c,j}}{\partial x^{c,j}_i}  \frac{\partial f^{c,j}_i}{\partial \omega_j}
+ \sum_i \frac{\partial g^{u,j}}{\partial x^{u,j}_i}  \frac{\partial f^{u,j}_i}{\partial \omega_j}  > 0. 
\end{align}
These conditions are, for example, satisfied by the demand dynamics considered in Section~\ref{Simulation}.
}
\end{remark}


In addition, we impose an assumption concerning the form of the cost functions in the OSLC problem~\eqref{Problem_To_Min}.

\begin{assumption} \label{assum5}
The cost functions $C_{j}$ and $C_{dj}$ are continuously differentiable and strictly convex. Moreover, the first derivative of $h_j^{-1}(z)$ is nonnegative for all $z \in \mathbb{R}$.
\end{assumption}

{\color{black}
\begin{remark}
The convexity and continuous differentiability of the cost functions given by Assumption~\ref{assum5} are sufficient to allow the use of the KKT conditions to prove the optimality result in Theorem~\ref{optthm}. We demonstrate in Appendix B how the condition of continuity in the derivatives may be relaxed by {\color{black}by means of subgradient methods.}
\end{remark}
}

We will explicitly state within each of our results which of these conditions are required.

\section{Main results} \label{Main_Results}

In this section we state our main results, with their {\color{black}proofs provided} in Appendix A. Our first result shows that the set of equilibria of the system~\eqref{sys1}--\eqref{sys2} for which the assumptions stated are satisfied
 is asymptotically attracting, while our second result demonstrates sufficient conditions for equilibrium points to be optimal for the OSLC problem~\eqref{Problem_To_Min}. Based on these results, we can guarantee convergence to optimality of all solutions starting in the vicinity of an equilibrium. Finally, we show that the inclusion of controllable demand in our model reduces steady state frequency deviation, thereby aiding in frequency control.

\begin{theorem}[\color{black}Stability] \label{convthm}
Suppose that Assumptions~\ref{assum1}--\ref{assum4} are all satisfied.
Then there exists an open neighbourhood $S$ of the equilibrium $(\eta^*, \omega^{G,*}, x^{M,*}, x^{c,*}, x^{u,*})$ such that whenever the initial conditions $(\eta(0), \omega^G(0),x^M(0), x^c(0),x^u(0)) \in S$, then the solutions of the system~\eqref{sys1}--\eqref{sys2}
converge to an equilibrium as defined in Definition~\ref{eqbrdef}.
\end{theorem}

\begin{remark}
{\color{black}It will be seen within the proof of Theorem~\ref{convthm} that $\omega, x^{M}, x^{c}, x^{u}$ converge to $\omega^*, x^{M,*}, x^{c,*}, x^{u,*}$ respectively.
The phase differences $\eta$ also converge to a fixed point, however, this can be different from~$\eta^*$.
}
\end{remark}

{\color{black}
\begin{theorem}[\color{black}Optimality] \label{optthm}
Suppose that Assumption~\ref{assum5} is satisfied. If the control dynamics in~\eqref{sys2p} and~\eqref{sys2dc} are chosen such that {\color{black}$k_{p^M_j}(-\omega^*_j) = [(C_j')^{\hspace{-0.5pt}-1}(-\omega^*_j)]^{p^{M,max}_j}_{p^{M,min}_j}$} and $k_{d^c_j}(-\omega^*_j) = [(C_{dj}')^{\hspace{-0.5pt}-1}(\omega^*_j)]^{d^{c,max}_j}_{d^{c,min}_j}$, then the values $p^{M,*}$, $d^{c,*}$, and $d^{u,*}$ are optimal for the OSLC problem~\eqref{Problem_To_Min}.
\end{theorem}
}

From Theorems~\ref{convthm} and~\ref{optthm}, we immediately deduce the following result guaranteeing convergence to optimality.

{\color{black}\begin{theorem}[\color{black}Convergence to optimality] \label{mrthm}
Consider equilibria of~\eqref{sys1}--\eqref{sys2} with respect to which Assumptions~\ref{assum1}--\ref{assum5} are all satisfied. If the control dynamics in~\eqref{sys2p} and~\eqref{sys2dc} are chosen such that
\begin{equation}
{\color{black}\begin{aligned}
&k_{p^M_j}(\bar{u}) = [ (C_j')^{\hspace{-0.5pt}-1}(\bar{u})]^{p^{M,max}_j}_{p^{M,min}_j} \\
&k_{d^c_j}(\bar{u}) = [(C_{dj}')^{\hspace{-0.5pt}-1}(-\bar{u})] ^{d^{c,max}_j}_{d^{c,min}_j}
\label{contspec}
\end{aligned}}
\end{equation}
hold for all $\bar{u} \in \mathbb{R}$, then there exists an open neighbourhood of initial conditions about any such equilibrium such that the solutions of~\eqref{sys1}--\eqref{sys2} are guaranteed to converge to a global minimum of the OSLC problem~\eqref{Problem_To_Min}.
\end{theorem}}

\begin{remark}
Theorem~\ref{mrthm} states that if the system~\eqref{sys1}--\eqref{sys2} starts sufficiently close to any of its equilibria with respect to which {\color{black}Assumptions \ref{assum1} - \ref{assum5}} are satisfied, then the system is guaranteed to converge to an equilibrium point which will be optimal for the desired OSLC problem~\eqref{Problem_To_Min}. The fact that $p^M$ and $d^c$ represent controllable quantities means that we are free to design the dynamics in~\eqref{sys2p} and~\eqref{sys2dc} in order that the conditions~\eqref{contspec} are satisfied. Thus, knowledge of the cost functions in the optimization problem we want to solve explicitly {\color{black}determines} classes of dynamics which are guaranteed to yield convergence to optimal solutions.
\end{remark}

\begin{theorem}[{\color{black}Reduction in steady state error}]
\label{Theorem3}
Suppose that Assumptions~\ref{assum1}--\ref{assum5} are all satisfied with respect to all equilibria of~\eqref{sys1}--\eqref{sys2}. If the control dynamics in~\eqref{sys2p} and~\eqref{sys2dc} are chosen such that~\eqref{contspec} hold for all $\bar{u} \in \mathbb{R}$, then the addition of controllable demand in primary control results in a drop in steady state frequency deviation from its nominal value.
\end{theorem}

\section{Discussion} \label{Discussion}
\hspace{-1.5mm}
We now discuss various examples of generation and load dynamics that can fit within our framework.

As a first example, consider the model in~\cite{P19}, which investigates a linearized version of the system~\eqref{sys1} coupled with the static nonlinearities  {\color{black}$d^c_j = (C_{dj}')^{\hspace{-0.5pt}-1}(\omega_j)$} {\color{black}for the controllable demand, and with uncontrollable loads of the form $d^u_j = D_j \omega_j$. The damping constants $D_j$ were assumed positive, the cost functions $C_{dj}$ were taken to be strictly convex, and the mechanical power injection $p^M$ was also assumed to be constant after a step change. It is easy to see that for such a system Assumptions~\ref{assum1}--\ref{assum5} are all satisfied. Hence, this model can be analyzed in the framework introduced {\color{black}above, thus  implying} optimality and stability of the equilibrium points.}

\begin{figure}[t]
\centering
\includegraphics[trim = 43mm 90mm 43mm 100mm,width=3.115in,clip=true]{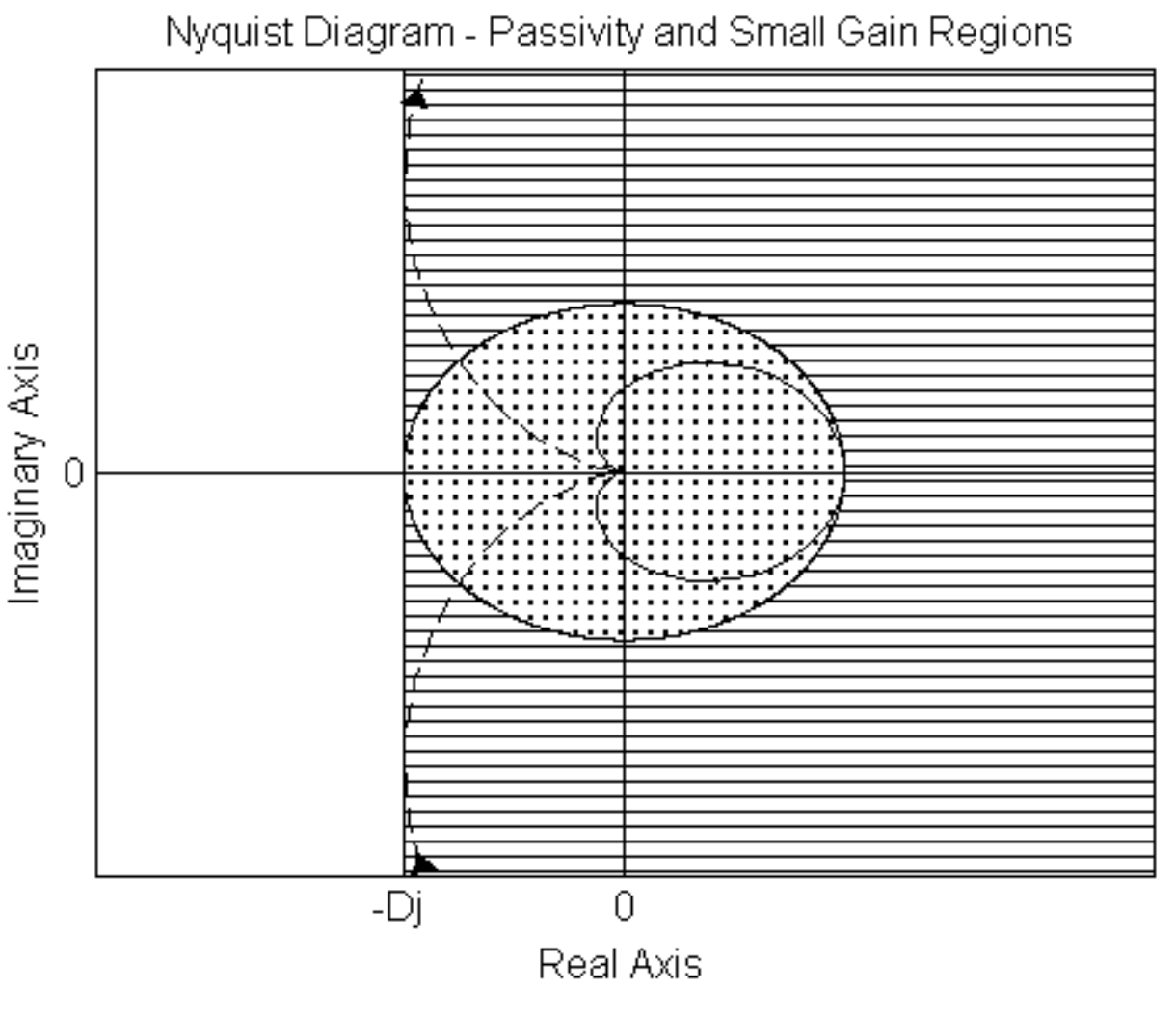}
\vspace{-2mm}
\caption{Nyquist plot for the transfer function relating {\color{black}$\hat{\tilde{p}}^M_j$} with {\color{black}$-\hat{\tilde{\omega}}_j$} when a linearization of~\eqref{9.1} about equilibrium is considered. A transfer function with a Nyquist plot within the circle satisfies the gain condition~\eqref{9.3}. Our approach allows the Nyquist plot to extend within the entire striped region.} \vspace{-5mm}
\label{Nyquist}
\end{figure}

The present framework can also include systems in which the generated powers satisfy any first-order dynamics, since such schemes {\color{black}are passive about their equilibria} for arbitrary gains. For {\color{black}higher-order} schemes, however, the dynamics for~$p^M$ are not necessarily passive, so some additional conditions are needed to ensure stability. {\color{black} Second-order generation dynamics are often considered in {\color{black}the} literature to model {\color{black}turbine-governor} dynamics {\color{black}(e.g.~\cite{Bergen_Vittal})}. {\color{black}These can be described by the system}}
\begin{equation} \label{9.1}
\begin{aligned}
&\dot{\alpha}_j = -\frac{1}{\tau _{g,j}} \alpha_j + \frac{1}{\tau _{g,j}} p_j^c, \\
&\dot{p}^M_j = -\frac{1}{\tau _{b,j}} p^M_j + \frac{1}{\tau _{b,j}} \alpha_j,
\end{aligned} \hspace{2em} j \in G
\end{equation}
 where $\alpha_j$ is the valve position of the turbine, the constants $\tau _{g,j}$ and $\tau _{b,j}$ represent lags in the dynamics of the governor and turbine respectively, and $p^c_j$ is a static function of frequency, {\color{black} corresponding to droop control}.
{\color{black}
{\color{black} Consider the case where} generator damping and {\color{black}uncontrollable} {\color{black}frequency-dependent} loads are {\color{black}modeled} by}
\begin{equation} \label{Uncontrollable_loads}
{\color{black}d^u_j = D_j \omega_j, \forall j \in N,}
\end{equation}
{\color{black} and there is no controllable demand.}  {\color{black} In~\cite{P34 - Low}, the gain condition
\begin{equation}
{\color{black}|p_j^c(\omega_j) - p_j^c(\omega_j^*)| \leq K_j |\omega_j - \omega_j^*|, \; j\in G}
\label{9.3}
\end{equation}
 with $K_j < D_j$ was imposed to ensure stability of {\color{black} the network when the generations dynamics are given by~\eqref{9.1}--\eqref{Uncontrollable_loads}}.}

{\color{black} As shown in {\color{black}Corollary}~\ref{Small_Gain_on_Low} below, under~\eqref{9.3}}, the overall system relating {\color{black}$-\omega_j$ with\footnote{\color{black}Note that {\color{black}this example could also include passive controllable demand $d^c_j(-\omega_j)$,
since} showing input strict passivity about equilibrium of the 
system with input $-\omega_j$ and output
$s^G_j = p^M_j - d^u_j$ is sufficient to ensure also that the system with the same input and output $s^G_j = p^M_j - d^u_j- d^c_j$ is input strictly passive.}
$s^G_j = p^M_j - d^u_j$ becomes input strictly passive about any equilibrium point. {\color{black}This follows from a more general result which we now state describing the connection between the $\mathcal{L}_2$-gain of general generation dynamics and the passivity of the supply dynamics. The proofs of Proposition~\ref{Small_Gain_on_Low} and Corollary~\ref{Small_Gain_on_Low} can be found in Appendix~A. }

{\color{black}
\begin{proposition}
Let equation~\eqref{Uncontrollable_loads} hold and consider any generation dynamics from {\color{black}$-\omega_j$} to $p^M_j$ {\color{black}of the form \eqref{sys2p}}. {\color{black} Given any equilibrium,} if the $\mathcal{L}_2$-gain from $(\omega_j - \omega_j^*)$ to $(p^M_j - p^{M,*}_j)$ is strictly less than $D_j$, then the system with input $-\omega_j$ and output {\color{black}$s^G_j = p^M_j - d^u_j$} is input strictly passive {\color{black}about the equilibrium considered.}
\label{Small_Gain}
\end{proposition}
}
{\color{black}
It easy to show {\color{black}that,} for the {\color{black}dynamics~\eqref{9.1}--\eqref{Uncontrollable_loads}}, condition~\eqref{9.3} is sufficient to satisfy the {\color{black}$\mathcal{L}_2$-gain} condition in Proposition~\ref{Small_Gain}. Applying Proposition~\ref{Small_Gain} on the system {\color{black}\eqref{9.1}--\eqref{9.3}} results in the following {\color{black}corollary}.
}
{\color{black}
\begin{corollary}
Let equation~\eqref{Uncontrollable_loads} hold and consider the generation dynamics in~\eqref{9.1}. Then, for any equilibrium where
~\eqref{9.3} holds, the system with input $-\omega_j$ and output $s^G_j = p^M_j - d^u_j$ is input strictly passive about this equilibrium.
\label{Small_Gain_on_Low}
\end{corollary}
}

%

{\color{black}Our framework can allow us to include also more general dynamics and {\color{black}to} deduce asymptotic stability under {\color{black}weaker conditions}.}
{\color{black}To} see this, we consider a linearization of the system~\eqref{9.1} {\color{black}about equilibrium
{\color{black} and let $\tilde{q}$ denote the deviation of any quantity $q$ from its equilibrium value $q^{\ast}$.} {\color{black}Expressing}} $\tilde{p}^M_j$ in the Laplace domain gives {{\color{black}
$\hat{\tilde{p}}_j^M = \frac{1}{(\tau _{g,j}s + 1)(\tau _{b,j}s + 1)} \hat{\tilde{p}}_j^c$}.
Therefore,
\begin{align}
\hat{\tilde{s}}^G_j = \hat{\tilde{p}}^M_j - \hat{\tilde{d}}^u_j &= \frac{1}{(\tau _{g,j}s + 1)(\tau _{b,j}s + 1)} \hat{\tilde{p}}_j^c + D_j (-\hat{\tilde{\omega}}_j) \nonumber \\
&=: {\color{black}H_j(s)} [-\hat{\tilde{\omega}}_j], \; j\in G, \label{9.5}
\end{align}
{\color{black}where $H_j(s)$ denotes the transfer} function relating $-\hat{\tilde{\omega}}_j$ and $\hat{\tilde{s}}^G_j$. Since the maximum gain of the transfer function from $\tilde{p}^c_j$ to $\tilde{p}^M_j$ is $1$ at $s = 0$, the condition in~\eqref{9.3} constrains the Nyquist diagram of $H_j$ to lie inside a ball with centre $(D_j,0)$ and radius $K_j < D_j$. This is contained strictly within the right half-plane, implying the required passivity condition in Assumption~\ref{assum2}. For instance, the Nyquist plot from input {\color{black}$-\hat{\tilde{\omega}}_j$} to output {\color{black}$\hat{\tilde{p}}^M_j$} can be as shown by the solid line in Fig.~\ref{Nyquist}. However, according to our analysis any dynamics for the command signal can be permitted provided that the supply dynamics in~\eqref{9.5} remain input-strictly passive. This can permit any frequency response within the striped region in Fig.~\ref{Nyquist}, for example allowing the larger Nyquist locus shown with a dashed line.
\color{black}{
In fact, under the {\color{black}reasonable} assumption that $p^c_j$ has the same sign as~$-\omega_j$ {(\color{black}i.e. negative feedback is used)}, 
{\color{black}it can easily be} verified that the transfer function from $\tilde{p}_j^c$ to $\tilde{p}_j^M$ {\color{black}given by $T_j(s)=\frac{1}{(\tau _{g,j}s + 1)(\tau _{b,j}s + 1)}$ }} has a minimum real value
\begin{equation} \label{mintf}
\Re\left({\color{black}T}_j(j\omega_{j,MAX}) \right)
= \frac{-\tau_{g,j}\tau_{b,j}}{(\tau_{g,j}+\tau_{b,j})^2 + 2(\tau_{g,j} + \tau_{b,j})\sqrt{\tau_{g,j}\tau_{b,j}}}
\end{equation}
{\color{black}at frequency $\omega_{j,MAX} = \sqrt{\frac{(\tau_{g,j} + \tau_{b,j}) + \sqrt{\tau_{b,j}\tau_{g,j}}}{(\tau_{b,j}\tau_{g,j})^{3/2}}}$. Thus, the required passivity property will be maintained provided $K_j$ multiplied by the quantity in~\eqref{mintf} is {\color{black} strictly} greater than $-D_j$. The maximum values of $K_j$ for which this is satisfied are shown in Fig.~\ref{Figure_with_K}, where $a$ is defined as the ratio $\tfrac{\tau_{b,j}}{\tau_{g,j}}$.} We see that the maximum allowable value for $K_j$ is always {\color{black}at least} $8D_j$ (obtained at $a = 1$) and tends to infinity as $a \rightarrow 0$ (which corresponds to a first order system). {\color{black}This shows that the stability guarantees can be preserved under significantly larger gains $K_j$ than those indicated in Corollary~\ref{Small_Gain_on_Low}. Therefore,} our approach allows for a less conservative stability condition for equilibrium points where a linearization is feasible, while also allowing to consider a wider class of generation dynamics.
}
{\color{black}It is worth noting, {\color{black}however}, that the use of the more conservative $\mathcal{L}_2$-gain condition described in Proposition~\ref{Small_Gain} {\color{black}would generally be expected to yield {\color{black}better robustness properties.} 
Such trade-offs between {\color{black}feedback gain} and stability margin {\color{black}need to be} taken into account in the design of control systems.}}

\begin{figure}[t!]
\centering
\includegraphics[trim = 35mm 90mm 43mm 95mm,width=3.3in,height = 2.3in,clip=true]{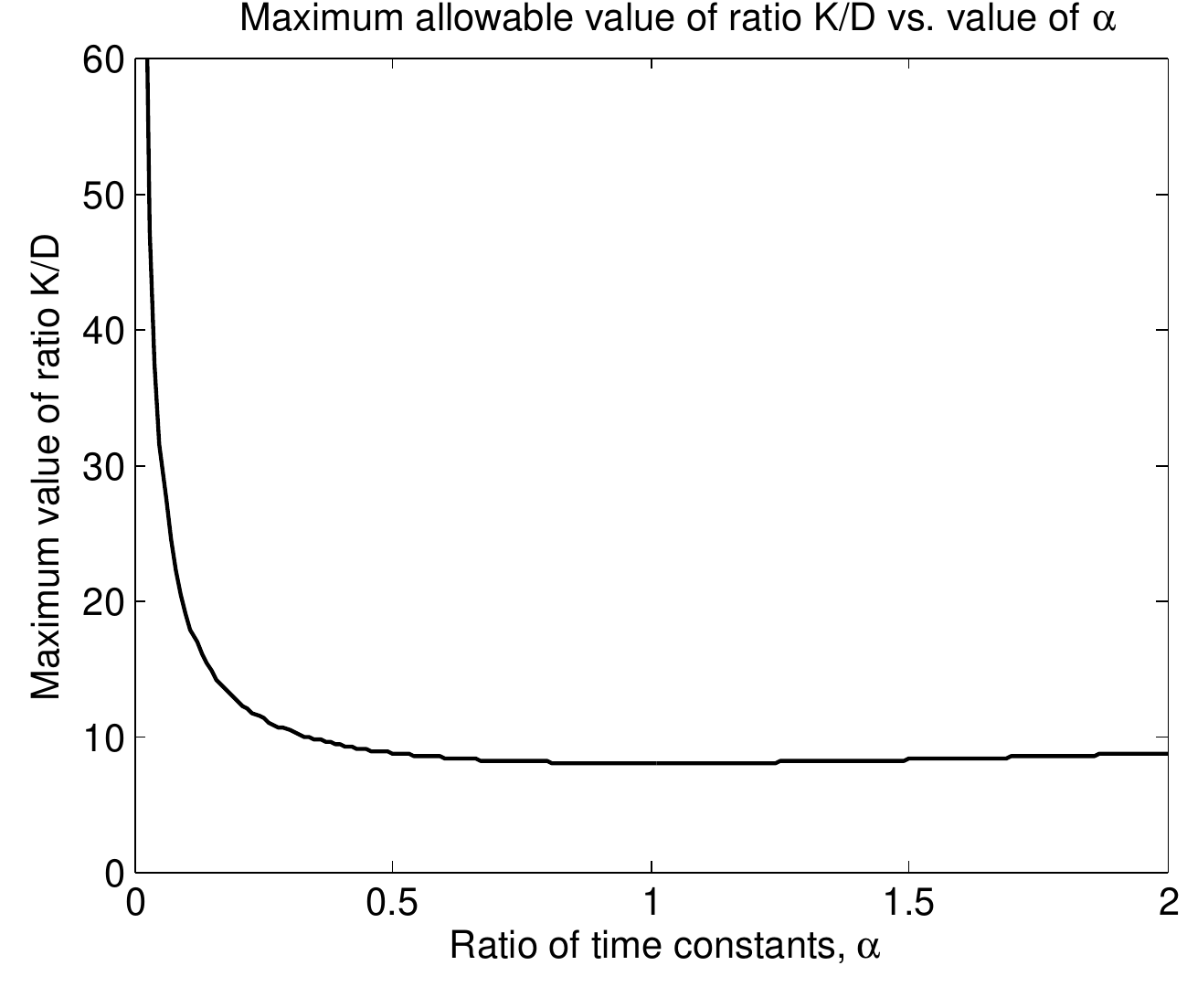}
\vspace{-2mm}
\caption{{\color{black} Maximum value of {\color{black}the} ratio $K_j/D_j$ under which the passivity property is maintained for the linearized dynamics of {\color{black}system~\eqref{9.1}--\eqref{9.3}} against the ratio of {\color{black}time constants $a=\tau_{b,j}/\tau_{g,j}$}. {\color{black}The figure} demonstrates that the maximum value for $K_j$ that ensures passivity is allowed to be much higher~than~$D_j$.}}
\label{Figure_with_K}
\vspace{-5mm}
\end{figure}

{\color{black}
\begin{remark}
{\color{black}
The assumption of small power angles differences $\eta$ is widely used in the literature.
Under such an assumption the stability results stated in the paper hold globally for the other system states if assumptions \ref{assum2}-\ref{assum4} hold globally\footnote{\color{black}In particular, system~\eqref{sys1}--\eqref{sys2} with the sinusoids in \eqref{sys1d} linearized about about $\eta^\ast$ (as in e.g.~\cite{P34 - Low}), would be globally asymptotically stable if in assumption \ref{assum2} the systems are globally passive about the equilibrium point considered, and in assumption \ref{assum3} the storage functions are radially unbounded with a minimum at the equilibrium point. For example, cases where the input strict passivity property is satisfied globally include the class of dynamics considered in Proposition~\ref{Small_Gain} (and Corollary \ref{Small_Gain_on_Low}), or the dynamics in \eqref{Static}--\eqref{Dynamic}.}.
It should be noted though that in general the system~\eqref{sys1}--\eqref{sys2} is not globally asymptotically stable 
due to the fact that the sinusoids in \eqref{sys1d} are only locally passive.
Finding the region of attraction is an interesting problem (associated with what is know as transient stability) which is not explicitly addressed within this paper.}
\end{remark}

{\color{black}
{\color{black} {\color{black} Another class of droop control schemes  
used in practice
incorporates a {deadband} which prevents unneeded adjustments for small variations in frequency about its nominal value. An example of this
is shown {in Fig.}~\ref{Deadband_system}. {The response profile of a static controllable load with deadband and saturation bounds is analogous.} For these systems, a minimum frequency deviation
$\omega_j^0$ is required to trigger a frequency-dependent deviation. }
{\color{black}The system then reaches its} physical limits at a higher frequency deviation $\omega_j^1$}.}

{\color{black}
While the stability of these systems can be shown by Theorem~\ref{convthm}, provided\footnote{\color{black}Although the nonlinearity depicted in Fig.~\ref{Deadband_system} is itself not input strictly passive about its equilibria, the presence of any frequency-dependent damping of the form~\eqref{Uncontrollable_loads} will be sufficient to ensure that Assumption~\ref{assum2} holds.} Assumptions~\ref{assum1}--\ref{assum4} hold, Theorem~\ref{optthm} cannot be applied, since the deadband makes the cost functions {\color{black}fail to satisfy continuously differentiability}, therefore not meeting the conditions of Assumption~\ref{assum5}.} {\color{black}Theorem~\ref{Optimality_switch} in Appendix B shows how Theorem~\ref{optthm} can be generalized to include non-continuously differentiable cost functions, thereby permitting the application of our results to analyze systems with {\color{black}{such~dynamics.}}}

\begin{figure}[t!]
\centering
\includegraphics[trim = 35mm 90mm 43mm 95mm,width=3.3in,height = 2.3in,clip=true]{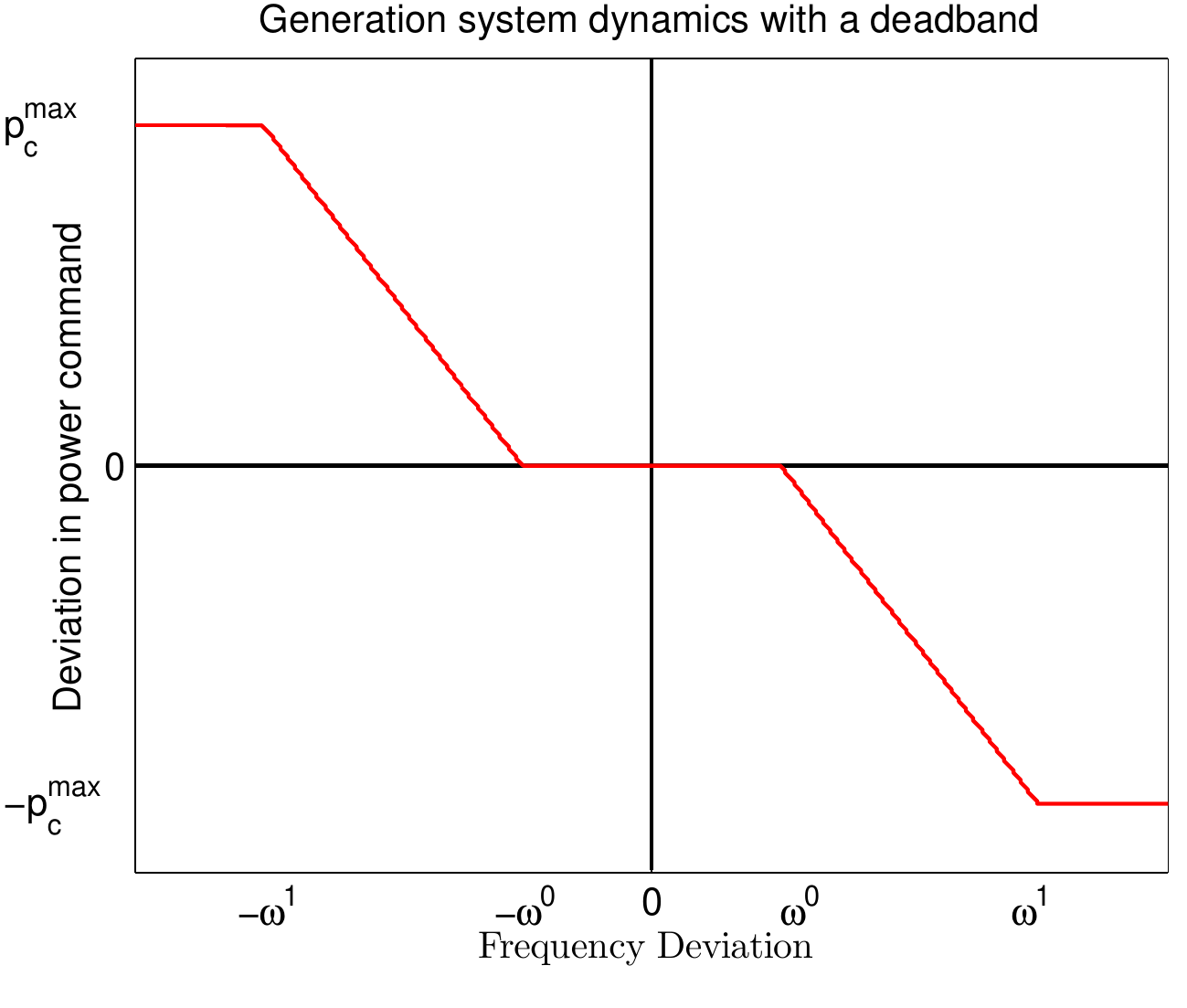}
\vspace{-2mm}
\caption{{\color{black} Power command {\color{black}signal} against frequency deviations for a system with deadband and saturation bounds.}}
\label{Deadband_system}
\vspace{-5mm}
\end{figure}

\section{Simulation on the IEEE 68 bus system} \label{Simulation}

In this section we illustrate our stability and optimality results through applications on the IEEE New York / New England 68 bus interconnection system~\cite{19-24}, simulated using the Power System Toolbox~\cite{19-25}. This model is more detailed and realistic than our {\color{black}analytical} model, including line resistances, a DC12 exciter model, power system stabilizer (PSS), and a subtransient reactance generator model. A similar model with no PSS was used for comparison\footnote{The details of the simulation models with or without PSS can be found in the Power System Toolbox data files data16m and data16em respectively.}.

The test system contains 52 load buses serving different types of loads including constant active and reactive loads. The overall system has a total real power of 16.41GW. For our simulation, we have added three loads on units 2, 9, and 17, each having a step increase of magnitude 1 p.u. (base 100MVA) at $t=1$ second. {\color{black} We allow controllable demand on buses 1-30. The disutility function for the aggregate load at each bus is $d^c_j$ is $C_{dj}(d^c_j) = \frac{1}{2}\alpha_j{d^{c}_j}^2$, where $\alpha_j$ are cost coefficients. The cost coefficients {\color{black}were selected such that} 
for a step change in demand, the power allocated between total generation and controllable demand would {\color{black}be roughly} equal, as suggested in~\cite{P50}. The selected values {\color{black}were} $\alpha_j = 5$ for load buses 1-10 and $\alpha_j = 10$ for the {\color{black}rest. The} loads were controlled every~10ms.}

Consider the static and dynamic control schemes given~by\footnote{Note that both~\eqref{Static} and~\eqref{Dynamic} are input strictly passive {\color{black}about the equilibria}, and both satisfy Assumption~\ref{assum4} (using the conditions in equations~\eqref{a4cond1} and~\eqref{a4cond2} respectively).}
\begin{equation}
d^c_j = (C_{dj}')^{\hspace{-0.5pt}-1} (\omega_j), \; j \in N,
\label{Static}
\end{equation}
\begin{equation}
\dot{d}^c_j = - (C_j'(d^c_j) - \omega _j), \; j \in N.
\label{Dynamic}
\end{equation}
We refer to the resulting dynamics as Static OSLC and Dynamic OSLC respectively. We investigate the behavior~in~the following six cases: (i)~no OSLC, no PSS, (ii) Static OSLC, no PSS, (iii) Dynamic OSLC, no PSS, (iv) no OSLC, with PSS, (v) Static OSLC, with PSS, (vi) Dynamic OSLC, with PSS. The frequency dynamics for bus 63 are shown in Fig.~\ref{Freq_Volt_ALL}.

\begin{figure}[t!]
\centering
\includegraphics[trim = 33mm 95mm 43mm 90mm,width=3.3in,height = 2.3in,clip=true]{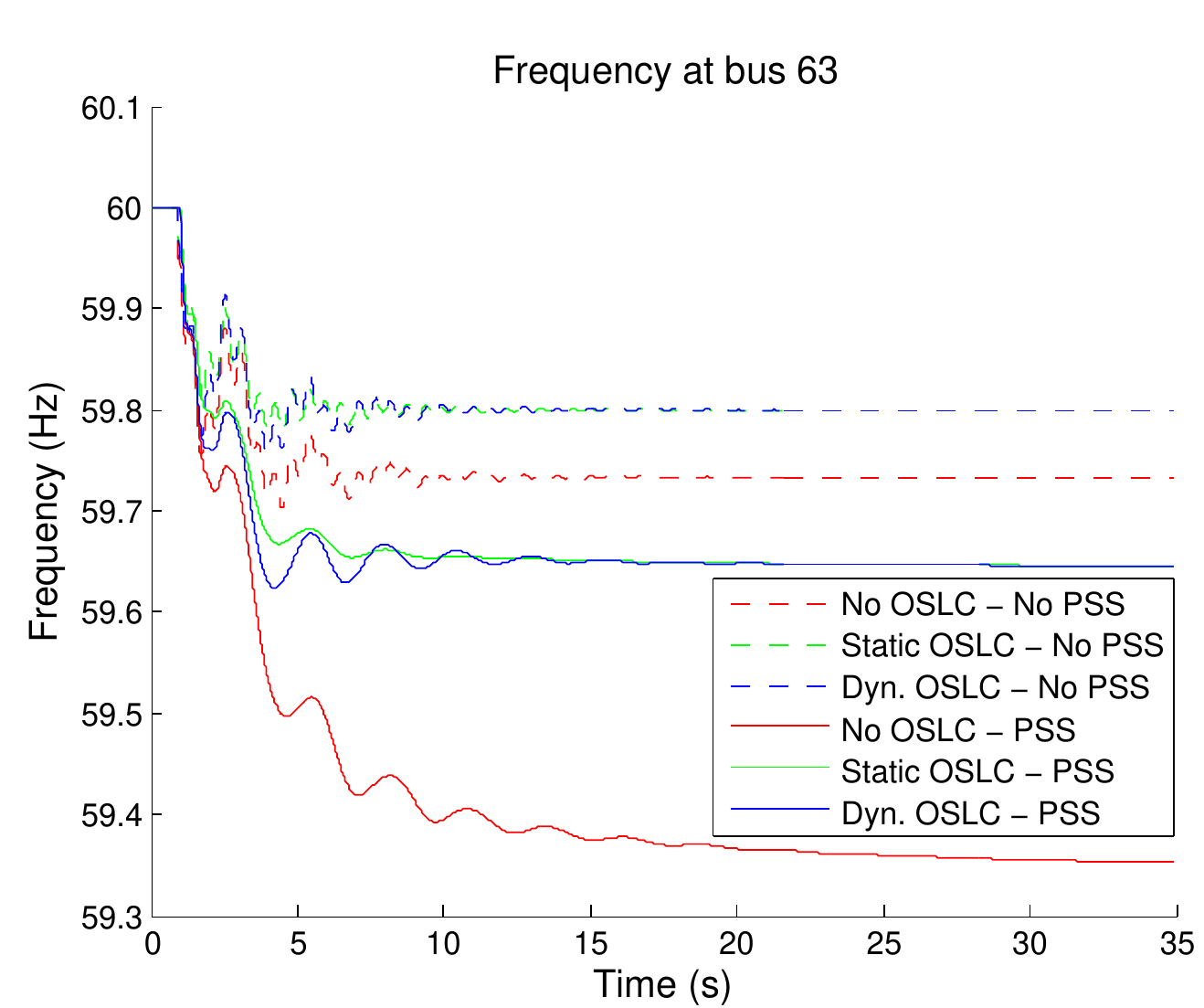}
\vspace{-2mm}
\caption{Frequency at bus 63 in six cases: (i)~no OSLC, no PSS, (ii) Static OSLC, no PSS, (iii) Dynamic OSLC, no PSS, (iv) no OSLC, with PSS, (v) Static OSLC, with PSS, (vi) Dynamic OSLC, with PSS.}
\label{Freq_Volt_ALL}
\vspace{-5mm}
\end{figure}

From Fig.~\ref{Freq_Volt_ALL}, we observe that whether or not PSS is used, the presence of OSLC results in a drop in steady state frequency deviations, as expected due to Theorem~\ref{Theorem3}. Furthermore, we can see that the overshoot is significantly less when OSLC is used.  The responses for Static and Dynamic OSLC have no significant differences and converge to the same exact value at steady state. However, Dynamic OSLC appears to give a larger overshoot than Static OSLC. {\color{black}{\color{black}In} all cases, the voltage deviation was less than 0.01 p.u., showing that the constant voltage assumption is {\color{black}reasonable}.}

{\color{black}In Fig.~\ref{Allocation_new} we also observe a higher power allocation at the load buses whose cost coefficients take the lower value $\alpha_j = 5$ than at those with $\alpha_j = 10$. This demonstrates that the power allocation among controllable loads depends upon the loads' respective cost coefficients of demand response. This behavior could be
beneficial if a 
{\color{black} prescribed allocation} were desirable, as then the load dynamics could be designed
{\color{black}such that} the
cost coefficients chosen yield the desired allocation.} {\color{black}Furthermore, as shown {\color{black}in} Fig.~\ref{Marginal_Costs}, the marginal costs {\color{black}at each controlled load} converge to the same value. {\color{black} This illustrates the optimality in the power allocation among the loads, as equality of marginal costs is the optimality condition for \eqref{Problem_To_Min} when the allocations do not saturate.}


\begin{figure}[t!]
\centering
\includegraphics[trim = 35mm 95mm 43mm 90mm,width=3.3in,height = 2.3in ,clip=true]{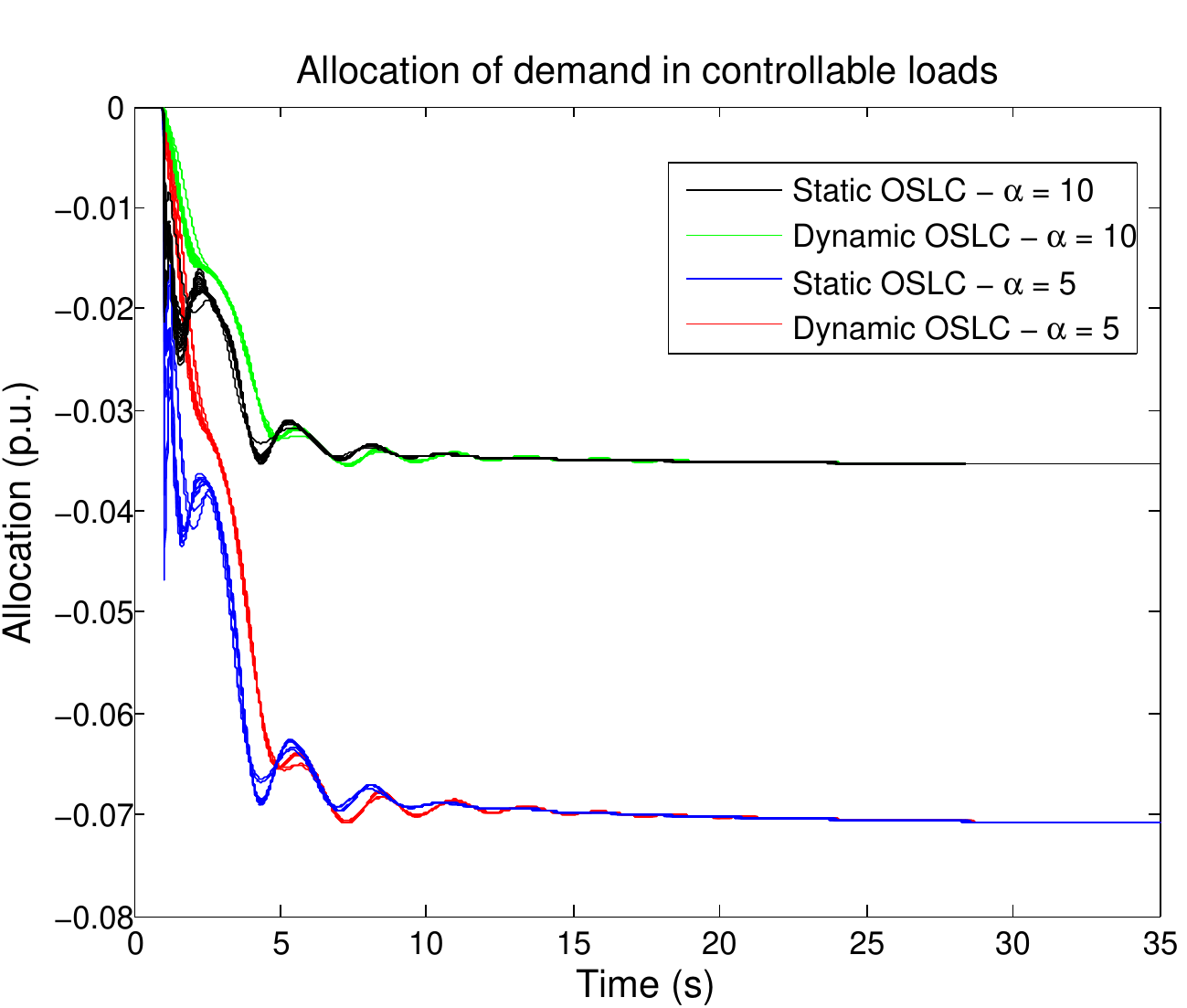}
\vspace{-2mm}
\caption{{\color{black}Power allocation among controllable loads {\color{black}with} non-equal cost coefficients for 2 cases: (i) Static OSLC (ii) Dynamic OSLC.}}
\label{Allocation_new}
\vspace{-5mm}
\end{figure}

\begin{figure}[t!]
\centering
\includegraphics[trim = 35mm 95mm 43mm 90mm,width=3.3in,height = 2.3in ,clip=true]{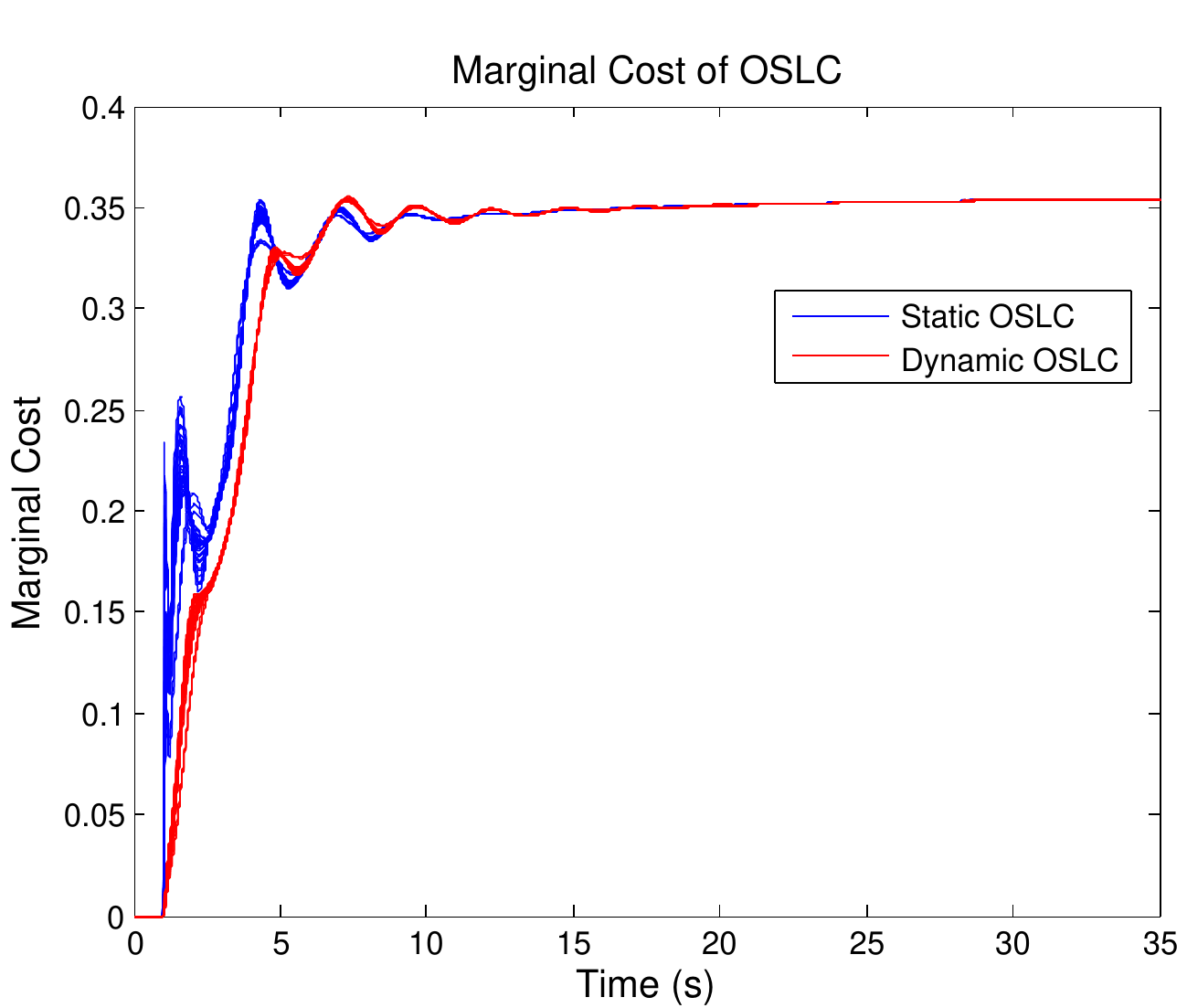}
\vspace{-2mm}
\caption{{\color{black}Marginal costs $C_{dj}'$ of controllable loads {\color{black}with} non-equal cost coefficients for 2 cases: (i) Static OSLC (ii) Dynamic OSLC.}}
\label{Marginal_Costs}
\vspace{-5mm}
\end{figure}


Finally, to investigate the system's robustness, we then introduced delays to account for the time between the arrival of the frequency signal and the response of the controllable demand. The simulation was repeated with 0.1 pu loads and a delay of 0.05 seconds. {\color{black}Furthermore, all cost {\color{black}coefficients} were set to $\alpha_j = 1$.} Dynamic OSLC was seen to offer improved robustness to the time-delay relative to Static OSLC, since the first converged both with and without PSS whereas the latter became unstable in both cases. This illustrates how appropriate {\color{black}higher order dynamics can have improved robustness properties}.
The simulation results for Dynamic OSLC are {\color{black}depicted in} Fig.~\ref{Freq_Volt_DEL}. }{\color{black}This} enhanced robustness to delays can be explained with the help of Fig.~\ref{Nyquist_Delay}.} {\color{black} The figure shows the Nyquist plots of the transfer functions  relating the increments from equilibrium of {\color{black}$d^c_j$ and $-\omega_j$}  for the systems~\eqref{Static} and~\eqref{Dynamic} respectively linearized about equilibrium,
with an input delay included and also both multiplied by a positive gain {\color{black}$K > D_j$}, where $D_j$ is the damping coefficient at bus $j$, as in \eqref{Uncontrollable_loads}}. 
{\color{black}The delayed Dynamic OSLC (dashed line)} 
{\color{black}maintains the passivity property of the bus dynamics
{\color{black}(since the Nyquist plot remains to right side of {\color{black}$-D_j$}, as was previously {\color{black}discussed in section \ref{Discussion}}).} {\color{black}On the other hand, the delayed Static OSLC (solid line) does not,} {\color{black}explaining why the latter might be expected to} become~unstable.}}

\begin{figure}[t]
\centering
\includegraphics[trim = 32mm 95mm 43mm 90mm,width=3.3in,height = 2.3in,clip=true]{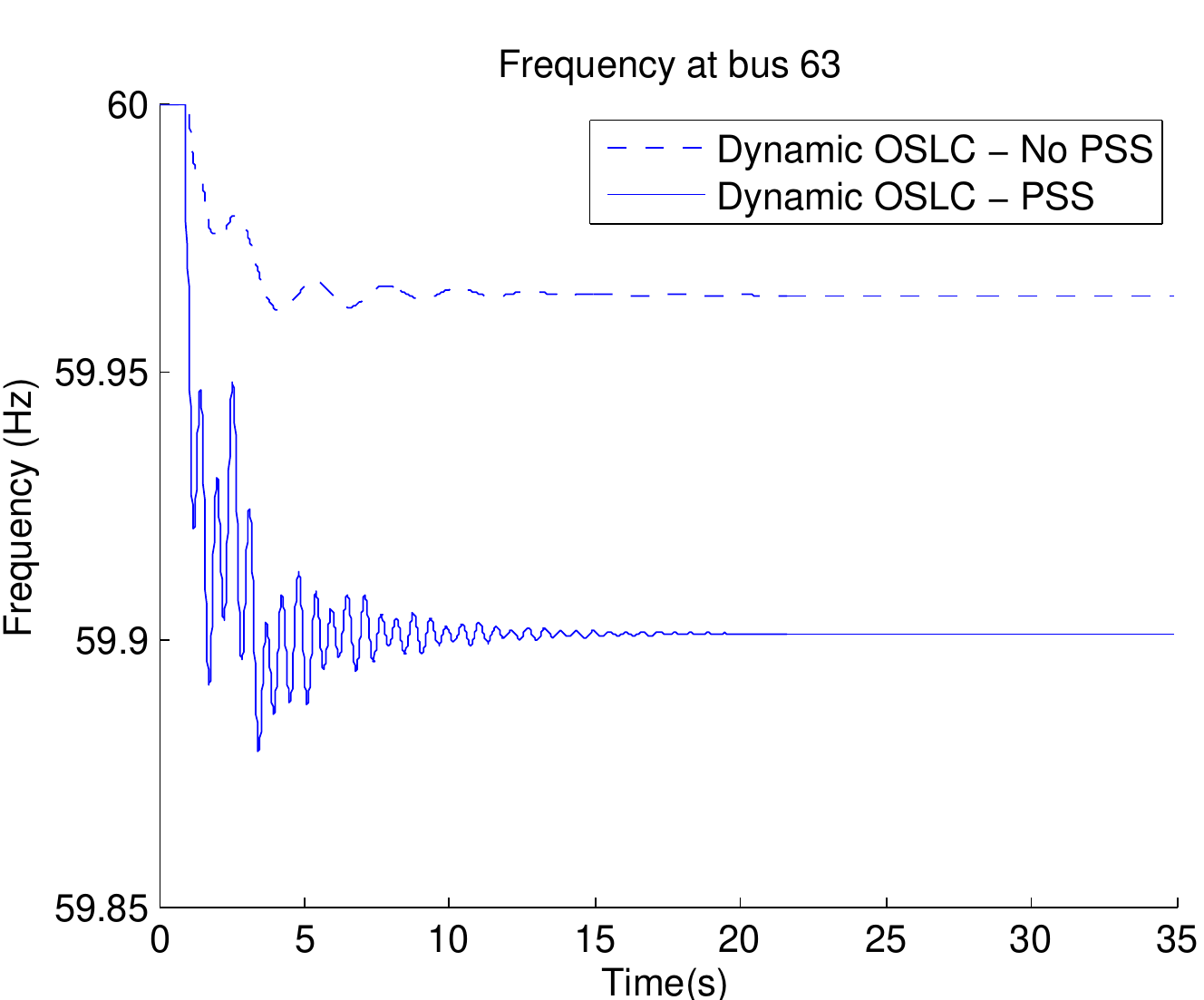}
\vspace{-2mm}
\caption{Frequency at bus 63 for Dynamic OSLC with time-delays in two cases: (i) no PSS (ii) with PSS.}
\label{Freq_Volt_DEL}
\vspace{-3mm}
\end{figure}

\begin{figure}[t!]
\centering
\includegraphics[trim = 43mm 93mm 43mm 100mm,width=3.3in,height = 2.4in,clip=true]{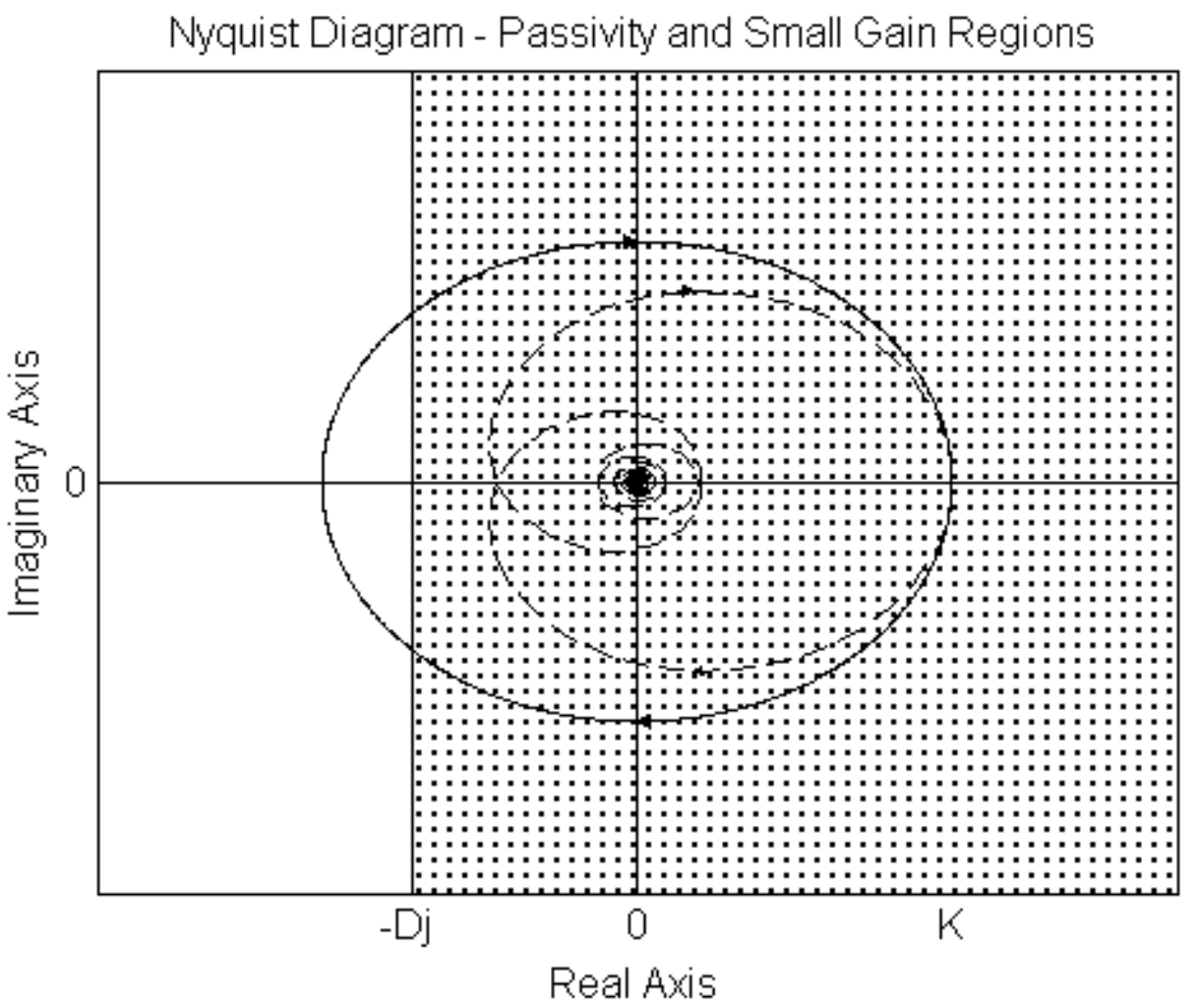}
\vspace{-2mm}
\caption{
{\color{black}Nyquist plot of the transfer functions relating the deviations from equilibrium of $d^c$ and $-\omega$, for the systems~\eqref{Static} and~\eqref{Dynamic} linearized about equilibrium, multiplied by a constant {\color{black}$K > D_j$} and with an input delay also included.}
{\color{black}The figure} shows that only the Dynamic delayed system (dashed line) remains on the right of {\color{black}$-D_j$}, {\color{black}maintaining} {\color{black}the passivity property of the bus dynamics when the damping coefficient is $D_j$ as in \eqref{Uncontrollable_loads}. The} Static system (solid line) {\color{black}extends} to the left of~{\color{black}{\color{black}$-D_j$}, {\color{black}hence the passivity property is lost}}.}
\label{Nyquist_Delay}
\vspace{-5mm}
\end{figure}
%
%

\section{Conclusion} \label{Discussion_Future_Work}

Within the paper, we have considered the problem of designing distributed generation and demand control schemes for primary frequency regulation in power networks such that asymptotic stability is guaranteed while ensuring optimality of power allocation. We have presented a network passivity framework which provides a systematic method to show stability over a broad class of generation and load dynamics. Furthermore, we have derived steady state conditions for the generation and controllable demand control schemes that ensure that the power generated/consumed is the solution to an appropriately constructed network optimization problem, thus allowing fairness in power allocation to be guaranteed. In addition, under some minor assumptions, we have shown that the inclusion within the model of controllable demand has a positive effect also on secondary control, decreasing the steady state deviation in frequency. A simulation on the IEEE 68 bus system verifies our results. Interesting potential extensions to our~analysis include incorporating voltage dynamics as well as more advanced models for the loads where their switching behavior is taken into account.
%
\linebreak

\section*{Appendix A \\ Proofs of our main results}\label{Appendix_A}

{\color{black}In this appendix we prove our main results, Theorems~\ref{convthm}--\ref{Theorem3}, and we also provide {\color{black}proofs} of Proposition~\ref{Small_Gain} {\color{black}and Corollary~\ref{Small_Gain_on_Low}}.}

\emph{Proof of Theorem~\ref{convthm}:}
We will use the dynamics in~\eqref{sys1} together with the passivity conditions in Assumption~\ref{assum2} to define a Lyapunov function for the system~\eqref{sys1}--\eqref{sys2}, similarly to \cite{P34 - Low}, \cite{trip}.

Firstly, we introduce $V_F (\omega^G) = \frac{1}{2}\sum_{j \in G} M_j (\omega_j - \omega^*_j)^2$. The time-derivative along the trajectories of~\eqref{sys1}--\eqref{sys2} is then
\begin{align*}
\dot{V}_F &= \sum_{j \in G} (\omega_j - \omega^*_j) \Bigg(-p^L_j + s_j^G - \sum_{k:j\rightarrow k} p_{jk} + \sum_{i:i\rightarrow j} p_{ij}\Bigg) \\
&\hspace{1em} + \sum_{j \in L} (\omega_j - \omega^*_j) \Bigg(-p^L_j + s_j^L - \sum_{k:j\rightarrow k} p_{jk} + \sum_{i:i\rightarrow j} p_{ij}\Bigg),
\end{align*}
by substituting~\eqref{sys1b} for $\dot{\omega}_j$ for $j \in G$ and adding the final term, which is equal to zero by~\eqref{sys1c}. Subtracting the product of $(\omega_j - \omega^*_j)$ with each term in~\eqref{eqbr2} and~\eqref{eqbr3}, this becomes
\begin{align}
&\dot{V}_F = \sum_{j \in G} (\omega_j - \omega^*_j) (s^G_j - s^{G,*}_j) + \sum_{j \in L} (\omega_j - \omega^*_j) (s^L_j - s^{L,*}_j) \nonumber \\
&\hspace{8em} + \sum_{(i,j) \in E} (p_{ij} - p^*_{ij}) (\omega_j - \omega_i), \label{VFdiff}
\end{align}
using in the final term the equilibrium condition~\eqref{eqbr1}.

Additionally, consider $V_P(\eta) = \sum_{(i,j) \in E} B_{ij} \int_{\eta^*_{ij}}^{\eta_{ij}} ( \sin \phi - \sin \eta^*_{ij} ) \, d\phi$. Using~\eqref{sys1a} and~\eqref{sys1d}, the time-derivative equals
\begin{align}
\dot{V}_P &= \sum_{(i,j) \in E} B_{ij} (\sin \eta_{ij} - \sin \eta^*_{ij}) (\omega_i - \omega_j) \nonumber \\
&= \sum_{(i,j) \in E} (p_{ij} - p^*_{ij}) (\omega_i - \omega_j). \label{VPdiff}
\end{align}

{\color{black}
Furthermore, from the passivity conditions in Assumption~\ref{assum2} and the definition of input strict passivity from Section~\ref{sec: Preliminaries}, it follows that the following hold:
\begin{enumerate}
\item For each $j \in G$, there exist open neighbourhoods $U^G_j$ of $\omega^*_j$ and $X^G_j$ of $(x^{M,j,*},x^{c,j,*},x^{u,j,*})$ and a continuously differentiable, positive semidefinite function\footnote{The time-derivates of $V^G_j$ and $V^L_j$ are evaluated using the dynamics~\eqref{sys2}.} $V^G_j(x^{M,j},x^{c,j},x^{u,j})$ such that
\begin{multline}\label{VG_derivative}
\dot{V}^G_j \le ((-\omega_j) - (-\omega^*_j))(s^G_j - s^{G,*}_j) \\
\hspace{8em} - \phi^G_j((-\omega_j) - (-\omega^*_j))
\end{multline}
for all $\omega_j \in U^G_j$ and all $(x^{M,j},x^{c,j},x^{u,j}) \in X^G_j$,
\item For each $j \in L$, there exist open neighbourhoods $U^L_j$ of $\omega^*_j$ and $X^L_j$ of $(x^{c,j,*},x^{u,j,*})$ and a continuously differentiable, positive semidefinite function $V^L_j(x^{c,j},x^{u,j})$ such that
\begin{multline}\label{VL_derivative}
\dot{V}^L_j \le ((-\omega_j) - (-\omega^*_j))(s^L_j - s^{L,*}_j) \\
\hspace{8em} - \phi^L_j((-\omega_j) - (-\omega^*_j))
\end{multline}
for all $\omega_j \in U^L_j$ and all $(x^{c,j}, x^{u,j}) \in X^L_j$,
\end{enumerate}
where $\phi^G_j$ and $\phi^L_j$ are positive definite functions.
}

Based on the above, we define the function
\begin{align*}
&V(\eta, \omega^G, x^M, x^c, x^u) = V_F(\omega^G) + V_P(\eta) \\
&\hspace{3em} + \sum_{j \in G} V^G_j(x^{M,j},x^{c,j},x^{u,j}) + \sum_{j \in L} V^L_j(x^{c,j},x^{u,j}),
\end{align*}
which we aim to use in Lasalle's theorem. By~\eqref{VFdiff} and~\eqref{VPdiff},
\begin{align*}
&\dot{V} = \sum_{j \in G} \left[(\omega_j - \omega^*_j) (s^G_j - s^{G,*}_j) + \dot{V}^G_j \right] \\
&\hspace{4em} + \sum_{j \in L} \left[(\omega_j - \omega^*_j) (s^L_j - s^{L,*}_j) + \dot{V}^L_j \right].
\end{align*}

Using~\eqref{VG_derivative} and~\eqref{VL_derivative} {\color{black}it therefore holds} that

\begin{align}
\dot{V} &\le -\sum_{j \in G} \phi^G_j((-\omega_j) - (-\omega^*_j)) - \sum_{j \in L} \phi^L_j((-\omega_j) - (-\omega^*_j)) \nonumber \\
&\le 0 \label{vdotineq}
\end{align}
whenever $\omega_j \in U^G_j$ and $(x^{M,j},x^{c,j},x^{u,j}) \in X^G_j$ for all $j \in G$ and $\omega_j \in U^L_j$ and $(x^{c,j},x^{u,j}) \in X^L_j$ for all $j \in L$.

Clearly $V_F$ has a strict global minimum at $\omega^{G,*}$, and $V^G_j$, $V^L_j$ have strict local minima at $(x^{M,j,*}, x^{c,j,*}, x^{u,j,*})$, $(x^{c,j,*}, x^{u,j,*})$ by Assumption~\ref{assum3}. Additionally, note that Assumption~\ref{assum1} guarantees the existence of some neighbourhood of each $\eta^*_{ij}$ on which the respective integrand in the definition of $V_P$ is increasing. Since the integrand is zero at the lower limit of the integration, $\eta^*_{ij}$, this immediately implies that $V_P$ has a strict local minimum at $\eta^*$. Thus, $V$ has a strict local minimum at the point $Q^* := (\eta^*, \omega^{G,*}, x^{M,*}, x^{c,*},  x^{u,*})$. We now recall Assumption~\ref{assum4}. This tells us that, provided $(\eta, \omega^G, x^M, x^c, x^u) \in T$, $\omega^L$ can be uniquely determined from these quantities. Therefore, the states of the differential equation system~\eqref{sys1}--\eqref{sys2} within the region $T$ can be expressed as $(\eta, \omega^G, x^M, x^c, x^u)$. Let us now choose a neighbourhood in the coordinates $(\eta, \omega^G, x^M, x^c, x^u)$ about $Q^*$ on which the following all hold:
\begin{enumerate}
\item $Q^*$ is a strict minimum of $V$,
\item $(\eta, \omega^G, x^M, x^c, x^u) \in T$,
\item $\omega_j \in U^G_j$ and $(x^{M,j},x^{c,j},x^{u,j}) \in X^G_j$ for all $j \in G$ and $\omega_j \in U^L_j$ and $(x^{c,j},x^{u,j}) \in X^L_j$ for all $j \in L$ \footnote{This is possible because the requirement $\omega_j \in U^L_j$ for all $j \in L$ corresponds, by Assumption~\ref{assum4} and the continuity of the equations in~\eqref{sys1}--\eqref{sys2}, to requiring the states $(\eta, \omega^G, x^M, x^c, x^u)$ to lie in some open neighbourhood about $Q^*$.},
\item $x^{M,j}$, $x^{c,j}$, and $x^{u,j}$ all lie within their respective neighbourhoods $X_0$ as defined in Section~\ref{sec: Preliminaries}.
\end{enumerate}
Recalling now~\eqref{vdotineq}, we thus see that within this neighbourhood, $V$ is a nonincreasing function of all the system states and has a strict local minimum at $Q^*$. Consequently, the connected component of the level set $\{(\eta, \omega^G, x^M, x^c, x^u) \colon V \le \epsilon\}$ containing $Q^*$ is guaranteed to be both compact and positively invariant with respect to the system~\eqref{sys1}--\eqref{sys2} for all sufficiently small $\epsilon > 0$. Therefore, there exists a compact positively invariant set $\Xi$ for~\eqref{sys1}--\eqref{sys2} containing~$Q^*$.

Lasalle's Invariance Principle can now be applied with the function $V$ on the compact positively invariant set~$\Xi$. This guarantees that all solutions of~\eqref{sys1}--\eqref{sys2} with initial conditions $(\eta(0), \omega^G(0), x^M(0), x^c(0) ,x^u(0)) \in \Xi$ converge to the largest invariant set within $\Xi \, \cap \, \{(\eta, \omega^G, x^M, x^c, x^u) \colon \dot{V} = 0\}$. We now consider this invariant set. If $\dot{V} = 0$ holds at a point within $\Xi$, then~\eqref{vdotineq} holds with equality, whence by Assumption~\ref{assum2} we must have $\omega = \omega^*$. Moreover, on any invariant set on which $\omega = \omega^*$, the system equations~\eqref{sys1} apply and give precisely the equilibrium conditions~\eqref{eqbr1},~\eqref{eqbr2},~\eqref{eqbr3}, and~\eqref{eqbr4}. {\color{black}Furthermore, if $\dot{V} = 0$, it follows from~\eqref{VFdiff},~\eqref{VPdiff},~\eqref{VG_derivative}, and~\eqref{VL_derivative} that all $\dot{V}^G_j = 0$ and $\dot{V}^L_j = 0$. But $\omega = \omega^*$ implies by the definitions in Section~\ref{sec: Preliminaries} that $(x^M, x^c, x^u)$ converge to the point $(x^{M,*}, x^{c,*}, x^{u,*})$, at which Assumption~\ref{assum3} states that $V^G_j$ and $V^L_j$ take strict local minima. Therefore the values of $V^G_j$ and $V^L_j$ must decrease along all nontrivial trajectories within the invariant set, {\color{black}contradicting $\dot{V}^G_j = 0$ and $\dot{V}^L_j = 0$}. Consequently, at all points of any invariant set within $\Xi \, \cap \, \{(\eta, \omega^G, x^M, x^c, x^u) \colon \dot{V} = 0\}$, we must also have $(x^M, x^c, x^u) = (x^{M,*}, x^{c,*}, x^{u,*})$. Thus, the remaining equilibrium conditions~\eqref{eqbr8},~\eqref{eqbr9},~\eqref{eqbr10},~\eqref{eqbr5},~\eqref{eqbr6}, and~\eqref{eqbr7} are also satisfied.} Therefore, we conclude by Lasalle's Invariance Principle that all solutions of~\eqref{sys1}--\eqref{sys2} with initial conditions $(\eta(0), \omega^G(0), x^M(0), x^c(0) ,x^u(0)) \in \Xi$ converge to the set of equilibrium points as defined in Definition~\ref{eqbrdef}. Finally, choosing for $S$ any open neighbourhood of $Q^*$ within $\Xi$ completes the proof. \hfill\IEEEQED

\emph{Proof of Theorem~\ref{optthm}:}
{\color{black}
Due to Assumption~\ref{assum5}, $C_j'$ and $C_{dj}'$ are strictly increasing and hence invertible. Therefore
\begin{equation}
{\color{black}\begin{aligned}
&k_{p^M_j}(-\omega_j^*) = [ (C_j')^{\hspace{-0.5pt}-1}(-\omega_j^*)]^{p^{M,max}_j}_{p^{M,min}_j} \\
&k_{d^c_j}(-\omega_j^*) = [(C_{dj}')^{\hspace{-0.5pt}-1}(\omega_j^*)] ^{d^{c,max}_j}_{d^{c,min}_j}
\label{pf2eq}
\end{aligned}}
\end{equation}
are well-defined. Furthermore, Assumption~\ref{assum5} also ensures that the OSLC {\color{black}problem}~\eqref{Problem_To_Min} is a convex optimization problem {\color{black}with a continuously differentiable cost function}. Thus, a point $(\bar{p}^M, \bar{d}^c, \bar{d}^u)$ is a global minimum for~\eqref{Problem_To_Min} if and only if it satisfies the KKT conditions~\cite{Boyd}
\begin{subequations} \label{kkt}
\begin{equation}
{\color{black}C_j'(\bar{p}^M_j) = -\nu - \lambda_j^+ + \lambda_j^-, \; j \in G,} \label{kkt1}
\end{equation}
\begin{equation}
{\color{black}C_{dj}'(\bar{d}^c_j) = \nu - \mu_j^+ + \mu_j^-, \; j \in N,} \label{kkt2}
\end{equation}
\begin{equation}
h_j^{-1}(\bar{d}^u_j) = \nu, \; j \in N, \label{kkt3}
\end{equation}
\begin{equation}
\sum\limits_{j\in  G} \bar{p}_j^M = \sum\limits_{j\in  N} (\bar{d}^c_j + \bar{d}^u_j +p_j^L), \label{kkt4}
\end{equation}
\begin{equation}
{\color{black}p^{M,min}_j \leq \bar{p}^M_j \leq p^{M,max}_j, \; j \in G,} \label{kkt5}
\end{equation}
\begin{equation}
{\color{black}d^{c,min}_j \leq \bar{d}^c_j \leq d^{c,max}_j, \; j \in N,}  \label{kkt6}
\end{equation}
\begin{equation}
{\color{black}\lambda_j^+ (p^M_j - p^{M,max}_j) = 0, \; \lambda_j^- (p^M_j - p^{M,min}_j) = 0, \; j \in G,} \label{kkt7}
\end{equation}
\begin{equation}
{\color{black}\mu^+ (d^c_j - d^{c,max}_j) = 0, \; \mu^- (d^c_j - d^{c,min}_j) = 0, \; j \in G,} \label{kkt8}
\end{equation}
\end{subequations}
{\color{black}for some constants $\nu \in \mathbb{R}$ and $\lambda_j^+, \lambda_j^-, \mu_j^+, \mu_j^- \ge 0$. We will now show that these conditions are satisfied by the equilibrium values $(\bar{p}^M, \bar{d}^c, \bar{d}^u) = (p^{M,*}, d^{c,*}, d^{u,*})$ defined by equations~\eqref{eqbr5},~\eqref{eqbr6}, and~\eqref{eqbr7}.

By Assumption~\ref{assum5}, $C_j'$ and $C_{dj}'$ are both strictly increasing, so we can uniquely define\footnote{\color{black}If any of these quantities do not exist, then we can just take $\omega_j^{M,max}, \omega_j^{c,min} := -\infty$ or $\omega_j^{M,min}, \omega_j^{c,max} := \infty$ as appropriate, and the following arguments still apply.} frequencies $\omega_j^{M,max} := -C_j'(p^{M,max}_j)$, $\omega_j^{M,min} := -C_j'(p^{M,min}_j)$, {\color{black}$\omega_j^{c,max} := C_{dj}'(d^{c,max}_j)$, and $\omega_j^{c,min} := C_{dj}'(d^{c,min}_j)$}. Letting $\omega^*_0$ denote the common value of all $\omega^*_j$ due to~\eqref{eqbr1}, we can, in terms of these quantities, define the nonnegative constants
\begin{equation} \label{lambdamudefs}
{\color{black}\begin{aligned}
&\lambda_j^+ := (\omega^{M,max}_j - \omega^*_0) \, \mathds{1}_{\{q \colon q \le \omega^{M,max}_j\}} (\omega_0^*), \\
&\lambda_j^- := (\omega^*_0 - \omega^{M,min}_j) \, \mathds{1}_{\{q \colon q \ge \omega^{M,min}_j\}} (\omega_0^*), \\
&\mu_j^+ := (\omega^*_0 - \omega^{c,max}_j) \, \mathds{1}_{\{q \colon q \ge \omega^{c,max}_j\}} (\omega_0^*), \\
&\mu_j^- := (\omega^{c,min}_j - \omega^*_0) \, \mathds{1}_{\{q \colon q \le \omega^{c,min}_j\}} (\omega_0^*).
\end{aligned}}
\end{equation}
{\color{black} Then, since $(C_j')^{\hspace{-0.5pt}-1}(-\omega_0^*) \ge p^{M,max}_j \Leftrightarrow \omega^*_0 \le \omega^{M,max}_j$, $(C_j')^{\hspace{-0.5pt}-1}(-\omega_0^*) \le p^{M,min}_j \Leftrightarrow \omega^*_0 \ge \omega^{M,min}_j$, $(C_{dj}')^{\hspace{-0.5pt}-1}(\omega_0^*) \ge d^{c,max}_j \Leftrightarrow \omega^*_0 \ge \omega^{c,max}_j$, and $(C_{dj}')^{\hspace{-0.5pt}-1}(\omega_0^*) \le d^{c,min}_j \Leftrightarrow \omega^*_0 \le \omega^{c,min}_j$, it follows by~\eqref{eqbr5},~\eqref{eqbr6}, and~\eqref{pf2eq} that the complementary slackness conditions~\eqref{kkt7} and~\eqref{kkt8} are satisfied with the choices~\eqref{lambdamudefs}.

Let us also define $\nu = \omega_0^*$. Then
\begin{align*}
(C_j')^{\hspace{-0.5pt}-1} ( -\nu - \lambda_j^+ + \lambda_j^- ) &= (C_j')^{\hspace{-0.5pt}-1} \Big( [-\omega_0^*]_{-\omega^{M,min}_j}^{-\omega^{M,max}_j} \Big) \\
&= [ (C_j')^{\hspace{-0.5pt}-1} (-\omega_0^*)]_{p^{M,min}_j}^{p^{M,max}_j} \\
&= p^{M,*}_j,
\end{align*}
by the above definitions and equations~\eqref{eqbr5} and~\eqref{pf2eq}. Therefore, the optimality condition~\eqref{kkt1} holds. Analogously,
\begin{align*}
(C_{dj}')^{\hspace{-0.5pt}-1} ( \nu - \mu^+ + \mu^- ) &= (C_{dj}')^{\hspace{-0.5pt}-1} \Big( [\omega_0^*]_{\omega^{c,min}_j}^{\omega^{c,max}_j} \Big) \\
&= [ (C_{dj}')^{\hspace{-0.5pt}-1} (\omega_0^*)]_{d^{c,min}_j}^{d^{c,max}_j} \\
&= d^{c,*}_j,
\end{align*}
by~\eqref{eqbr6} and~\eqref{pf2eq}, satisfying~\eqref{kkt2}. Additionally,~\eqref{kkt3} holds because $h_j(\nu) = d^u_j$ follows from~\eqref{eqbr7} and~\eqref{hdef}.}
%

Furthermore, summing the equilibrium conditions~\eqref{eqbr2} over all $j \in G$ and~\eqref{eqbr3} over all $j \in L$ shows that~\eqref{kkt4} {\color{black}is satisfied}. Finally, the saturation constrains in~\eqref{pf2eq} ensure that~\eqref{kkt5} and~\eqref{kkt6} are also satisfied.

Thus, we see that the values $(\bar{p}^M, \bar{d}^c, \bar{d}^u) = (p^{M,*}, d^{c,*}, d^{u,*})$ satisfy the KKT conditions~\eqref{kkt} with $\nu = \omega^*_0$. Therefore, the equilibrium values $p^{M,*}$, $d^{c,*}$, and~$d^{u,*}$ define a global minimum for the problem~\eqref{Problem_To_Min}. \hfill\IEEEQED
}
%

{\color{black}\emph{Proof of Theorem~\ref{mrthm}:}
If Assumptions~\ref{assum1}--\ref{assum5} all hold and~\eqref{contspec} is true, then all of the assumptions in both Theorems~\ref{convthm} and~\ref{optthm} are satisfied, and thus the result follows. \hfill \IEEEQED}

\emph{Proof of Theorem~\ref{Theorem3}:}
Recalling the proof of Theorem~\ref{optthm}, we know from~\eqref{kkt4} and the equalities~\eqref{pf2eq} that at any equilibrium of~\eqref{sys1}--\eqref{sys2} the power balance equation
\begin{equation}
 \sum \limits_{j\in G} (C_j')^{\hspace{-0.5pt}-1}(\omega^*_0) + \sum \limits_{j\in N} ( (C_{dj}')^{\hspace{-0.5pt}-1}(\omega^*_0) + h_{j}(\omega^*_0)) = - \sum \limits_{j\in N} p_j^L \label{pbal}
\end{equation}
is satisfied, where $\omega^*_0$ denotes the common steady state value of frequency due to~\eqref{eqbr1}. Now note that, because the nominal frequency defines an equilibrium frequency prior to the step change in load and all quantities in~\eqref{sys1} denote deviations from their respective values at this nominal equilibrium, the equalities~\eqref{hdef} and~\eqref{contspec} imply that each term on the left-hand side in~\eqref{pbal} must take the value zero at $\omega^*_0 = 0$. Furthermore, Assumption~\ref{assum5} implies that the terms $(C'_j)^{\hspace{-0.5pt}-1}(\omega^*_0)$ and $(C'_{dj})^{\hspace{-0.5pt}-1}(\omega^*_0)$ are all strictly increasing in $\omega^*_0$, while each term $h_j(\omega^*_0)$ is nondecreasing in $\omega^*_0$. Thus both the added term due to load control and the entire left-hand side in~\eqref{pbal} have the same sign {\color{black}as $\omega^*_0$} and are strictly increasing in $\omega^*_0$. It follows that the presence of this load control term results in a decrease in the value of $\omega^*_0$, the steady state frequency deviation from its nominal value. \hfill\IEEEQED

{\color{black}
\emph{Proof of Proposition~\ref{Small_Gain}:}
The $\mathcal{L}_2$-gain condition {\color{black}implies\footnote{\color{black} As in the main text, for a variable $x$ that depends on time we use the notation $\tilde{x}$ to denote deviations from equilibrium, i.e. $\tilde{x}:=x-x^\ast$, where $x^\ast$ is the value of $x$ at the equilibrium point mentioned in the proposition.}}
\begin{equation} \label{L2_gain_condition}
{\color{black}\sqrt{\int_0^{{t_1}} \! (\tilde{p}^{M}_j)^2 \, dt} \leq K_j \sqrt{\int_0^{{t_1}} \! \tilde{\omega}^2_j(t) \, dt.}}
\end{equation}
where $K_j < D_j$ {\color{black}and $t_1$ is any {\color{black}positive constant.}
Then, input strict passivity can be shown as follows.
\begin{subequations}\label{Pass_ineq_all}
\begin{align}
\int_0^{{t_1}} \! {\color{black}\tilde{p}^{M}_j(t) \tilde{\omega}_j(t)} \, dt &\leq \sqrt{\left(\int_0^{{t_1}} \! {\color{black}|\tilde{p}^{M}_j(t)| |\tilde{\omega}_j(t)|} \, dt \right)^2} \label{Pass_1} \\
&\hspace{-1em}{\color{black}\leq \sqrt{\int_0^{{t_1}} \! |\tilde{p}^{M}_j(t)|^2 \, dt \int_0^{{t_1}} \! |\tilde{\omega}_j(t)|^2 \, dt}} \label{Pass_2} \\
&\hspace{-1em}{\color{black}\leq \int_0^{{t_1}} \! K_j \tilde{\omega}_j(t)^2 \, dt
< \int_0^{{t_1}} \! D_j \tilde{\omega}_j(t)^2 \, dt
,} \label{Pass_3}
\end{align}
\end{subequations}
where inequality~\eqref{Pass_2} {\color{black}follows from the Cauchy–-Schwarz inequality} 
and~\eqref{Pass_3} from inequality~\eqref{L2_gain_condition} and {\color{black}$K_j < D_j$, $j \in G$}. Using~\eqref{Pass_ineq_all}, it is straightforward to show {\color{black}that}
\begin{multline}
\int_{0}^{t_1}D_j \tilde{\omega}_j(t)^2 dt - \int_{0}^{t_1} \tilde{p}^M_j(t) \tilde{\omega}_j(t) dt \\
{\color{black}=} \int_{0}^{t_1} \tilde{s}_j^G(t) (-\tilde{\omega}_j(t)) dt \geq (D_j - K_j) \int_{0}^{t_1}  \tilde{\omega}_j(t)^2 dt{\color{black}>0}
\label{Ap.n.4}
\end{multline}
holds for all {\color{black}$j \in G$. Inequality}~\eqref{Ap.n.4} implies input strict passivity {\color{black}of the system with output $\tilde{s}^G_j=\tilde{p}^M_j-\tilde{d}^u_j$ and input $-\tilde{\omega}_j$, about the equilibrium point considered, since \eqref{Ap.n.4} implies from~\cite[Lemma 1]{P28} the existence of} a positive definite storage function $V$ satisfying the input strict passivity condition {\color{black}in} Definition~\ref{Passivity_Definition}.~\hfill\IEEEQED

{\color{black}
\emph{Proof of {\color{black}Corollary}~\ref{Small_Gain_on_Low}:}
{\color{black}Let} ${\color{black}T}_j(s)$ be the transfer function relating {\color{black}$\hat{\tilde{p}}^c_j(s)$} and {\color{black}$\hat{\tilde{p}}^M_j(s)$} in~\eqref{9.1}, given by
\begin{equation*} \label{G_transfer_function}
{\color{black}T}_j(s) = \frac{1}{(\tau _{g,j}s + 1)(\tau _{b,j}s + 1)}, \; {\color{black}j \in G}.
\end{equation*}
It is easy to show that {\color{black}$\|{\color{black}T}_j(s)\|_{\infty} = 1$,} {\color{black}hence the system from $\tilde{p}^c_j$ to $\tilde{p}^M_j$ has $\mathcal{L}_2$-gain less than or equal to $1$ (e.g.~\cite[p.18]{Feedback_Control_Doyle}). Using also equation~\eqref{9.3} we thus have
\begin{equation} \label{Low_inequality}
\|\tilde{p}^M_j\|_2 \leq \|\tilde{p}_j^c\|_2 < D_j \|\tilde{\omega}_j\|_2.
\end{equation}
}
{\color{black} It follows from~\eqref{Low_inequality} that the $\mathcal{L}_2$-gain condition in Proposition~\ref{Small_Gain}} 
holds, 
{\color{black}therefore Proposition~\ref{Small_Gain} can be used
to deduce} input strict passivity of the system.
{\color{black}\hfill\IEEEQED}
}
{\color{black}
\begin{remark}
It can be verified that a possible choice of storage function for the input strictly passive system considered in Corollary~\ref{Small_Gain_on_Low} is
\begin{equation*}
{\color{black}V(p^M,\alpha) = \sum_{j\in G} \left( \frac{1}{2} \beta_j (p^M_j - p^{M,*}_j)^2 + \frac{1}{2} \gamma_j (\alpha_j - \alpha^*_j)^2 \right)}
\end{equation*}
with {\color{black}coefficients} $0 < \beta_j < \gamma_j \frac{\tau_{b,j}}{\tau_{g,j}}$ and $\gamma_j < \frac{2(D_j - K_j)}{K^2_j}   \tau_{g,j}$. 
Note that this has a strict minimum at the equilibrium point and hence also satisfies Assumption \ref{assum3}.
\end{remark}
}
\section*{Appendix B \\ Extension to non-$\mathcal{C}^1$ cost functions}\label{Appendix_B}

{\color{black}In this appendix we show how our main results can be extended to apply also to systems in which the cost functions within the OSLC problem~\eqref{Problem_To_Min} may not be continuously differentiable. We demonstrate this by formulating a relaxed version of Assumption~\ref{assum5}, {\color{black}which we use to}
to prove a generalization of our optimality result Theorem~\ref{optthm}. Since Theorem~\ref{convthm} is independent of Assumption~\ref{assum5} and so remains true also in this extended case, {\color{black}a corresponding generalization of }
Theorem~\ref{mrthm} will then follow immediately.

We weaken Assumption~\ref{assum5} to the following condition.

{\color{black}\begin{assumption} \label{assum5new}
The cost functions $C_{j}$ and $C_{dj}$ are continuous and strictly convex. Furthermore, they are continuously differentiable except on respective sets $\Lambda_j$ and $\Lambda_{dj}$ of isolated points. Moreover, the first derivative of $h_j^{-1}(z)$ is nonnegative for all $z \in \mathbb{R}$.
\end{assumption}

Assumption~\ref{assum5new} therefore relaxes the conditions imposed in the main paper by allowing the cost functions to be nondifferentiable at a discrete set of points. This permits classes of hybrid cost functions and more involved control dynamics. {\color{black}To} overcome the issue of nondifferentiability within the proof {\color{black}of the optimality result that follows}, we will consider {\color{black}the notions of subgradient and the subdifferential (e.g.~\cite{Rockafellar})}, defined as follows.
\begin{definition} \label{subdifdef}
Given a convex function $f : I \to \mathbb{R}$, a subgradient of $f$ at a point $x \in I \subseteq \mathbb{R}$ is any $v \in \mathbb{R}$ such that
$f(y) - f(x) \ge v (y - x)$  for all $y \in I$.
The set of all subgradients of $f$ at $x$ is called the subdifferential {\color{black}of $f$ at $x$ and is denoted by $\partial f (x)$.}
\end{definition}

{\color{black}As in the main text, the analysis needs to make use of the inverse of $C_j'$ and $C_{dj}'$. Since, however,  these derivatives are allowed to be discontinuous, their inverses are not well defined at points of discontinuity, and we hence need to define the following functions which can be seen as generalized inverses of $C_j'$ and $C_{dj}'$ respectively. In particular, we define}
\begin{equation} \label{geninverses}
\hspace{-0.1em}\begin{aligned}
D_j(x) &:= \hspace{-0.25em}\begin{cases}
(C_j')^{-1} (x), &\hspace{-0.7em}x \in C_j'(\mathbb{R} \setminus \Lambda_j), \\
\gamma, &\hspace{-0.7em}x \in [C_j'(\gamma-),C_j'(\gamma+)], \gamma \in \Lambda_j,
\end{cases} \hspace{-2em} \\
D_{dj}(x) &:= \hspace{-0.25em}\begin{cases}
(C_{dj}')^{-1} (x), &\hspace{-0.7em}x \in C_{dj}'(\mathbb{R} \setminus \Lambda_{dj}), \\
\gamma, &\hspace{-0.7em}x \in [C_{dj}'(\gamma-),C_{dj}'(\gamma+)], \gamma \in \Lambda_{dj},
\end{cases} \hspace{-2.5em}
\end{aligned}
\end{equation}
{\color{black} where the quantities $\gamma \pm$ respectively denote the limits $\lim_{\epsilon \downarrow 0} (\gamma \pm \epsilon)$.}
{\color{black} {\color{black}The functions $D_j$ and $D_{dj}$ defined in \eqref{geninverses} are therefore equal to {\color{black}$(C_{j}')^{-1}$ and $(C_{dj}')^{-1}$} respectively in the regime where $C_{j}'$ and $C_{dj}'$ are continuous, and otherwise remain constant.
They are thus continuous non decreasing functions analogous to {\color{black}$(C_{j}')^{-1}$ and $(C_{dj}')^{-1}$} and are well defined at all points.}
}

\begin{theorem} \label{Optimality_switch}
Suppose that Assumption~\ref{assum5new} is satisfied. If the control dynamics in~\eqref{sys2p} and~\eqref{sys2dc} are chosen such that
\begin{equation}
\begin{aligned}
&k_{p^M_j}(-\omega_j^*) = [ D_j(-\omega_j^*)]^{p^{M,max}_j}_{p^{M,min}_j}, \\
&k_{d^c_j}(-\omega_j^*) = [D_{dj}(\omega_j^*)] ^{d^{c,max}_j}_{d^{c,min}_j},
\label{pf5eq}
\end{aligned}
\end{equation}
then the values $p^{M,*}$, $d^{c,*}$, and $d^{u,*}$ are optimal for the OSLC problem~\eqref{Problem_To_Min}.
\end{theorem}}
}

{\color{black}
\emph{Proof of Lemma~\ref{Optimality_switch}:}
{\color{black}The optimality proof for non-zero values of generation or controllable demand is similar to that of Theorem~\ref{optthm}. {\color{black}We will make} use of {\color{black}subgradient techniques}~\cite[Section 23]{Rockafellar} {\color{black}to show that} the discontinuity in the cost function's derivative does not affect the optimality result.}

Firstly we note that strict convexity implies both that $C_j'$ and $C_{dj}'$ are strictly increasing on $\mathbb{R} \setminus \Lambda_j$ and $\mathbb{R} \setminus \Lambda_{dj}$ respectively, and that their jumps on $\Lambda_j$ and $\Lambda_{dj}$ are nonnegative. Therefore, $C_j'$ and $C_{dj}'$ are invertible on $\mathbb{R} \setminus \Lambda_j$ and $\mathbb{R} \setminus \Lambda_{dj}$, and {\color{black} $C_j'(\gamma-) < C_j'(\gamma+)$ for all $\gamma \in \Lambda_j$ and $C_{dj}'(\gamma-) < C_{dj}'(\gamma+)$ for all} {\color{black}$\gamma~\in~\Lambda_{dj}$}. Moreover, continuous differentiability on the sets $\mathbb{R} \setminus \Lambda_j$ and $\mathbb{R} \setminus \Lambda_{dj}$ ensures that the relevant limits $\gamma-$ and $\gamma+$ here all exist. Consequently, the functions in~\eqref{geninverses} and hence the controls in~\eqref{pf5eq} are well-defined.

Now, since OSLC~\eqref{Problem_To_Min} is a convex optimization problem, a point $(\bar{p}^M,\bar{d}^c,\bar{d}^u)$ is globally optimal if and only if it satisfies the subgradient KKT conditions. As only $C_j$ and {\color{black} $C_{dj}$ are} non-continuously differentiable, these are identical to those used in the proof of Theorem~\ref{optthm}, with the exception that~\eqref{kkt1} and~\eqref{kkt2} are replaced by the subdifferential optimality conditions
\begin{subequations}
\begin{equation}
-\nu - \lambda_j^+ + \lambda_j^- \in \partial C_j(\bar{p}^M_j), \; j \in G, \label{sdkkt1}
\end{equation}
\begin{equation}
\nu - \mu_j^+ + \mu_j^- \in \partial C_{dj}(\bar{d}^c_j), \; j \in N. \label{sdkkt2}
\end{equation}
\end{subequations}
We will thus seek to show that the arguments for demonstrating~\eqref{kkt3}--\eqref{kkt8} easily carry over to the present setting, and then prove that~\eqref{sdkkt1} and~\eqref{sdkkt2} are also satisfied by the equilibrium quantities, thereby guaranteeing optimality.

Due to the discontinuities in the cost function derivatives, we slightly redefine the quantities introduced in the proof of Theorem~\ref{optthm} as
\begin{equation} \label{omegaminmaxdefs}
\begin{aligned}
&\omega^{M,max}_j := \sup \{ \omega_j \colon D_j (-\omega_j) \ge p^{M,max}_j \},  \\
&\omega^{M,min}_j := \inf \{ \omega_j \colon D_j (-\omega_j) \le p^{M,min}_j \},  \\
&\omega^{c,max}_j := \inf \{ \omega_j \colon D_{dj} (\omega_j) \ge d^{c,max}_j \},  \\
&\omega^{c,min}_j := \sup \{ \omega_j \colon D_{dj} (\omega_j) \le d^{c,max}_j \}.
\end{aligned}
\end{equation}
It is easy to see that these quantities reduce to those considered in the proof of Theorem~\ref{optthm} in the case of continuously differentiable cost functions.

The definitions in~\eqref{omegaminmaxdefs} imply that $D_j(-\omega_0^*) \ge p^{M,max}_j \Leftrightarrow \omega^*_0 \le \omega^{M,max}_j$, $D_j(-\omega_0^*) \le p^{M,min}_j \Leftrightarrow \omega^*_0 \ge \omega^{M,min}_j$, $D_{dj}(\omega_0^*) \ge d^{c,max}_j \Leftrightarrow \omega^*_0 \ge \omega^{c,max}_j$, and $D_{dj}(\omega_0^*) \le d^{c,min}_j \Leftrightarrow \omega^*_0 \le \omega^{c,min}_j$. Therefore, the exact arguments in the proof of Theorem~\ref{optthm} also apply here and guarantee that~\eqref{kkt3}--\eqref{kkt8} all hold with the choices~\eqref{lambdamudefs}.

Finally, let $\nu = \omega_0^*$. Then
\begin{subequations}
\begin{align}
D_j ( -\nu - \lambda_j^+ + \lambda_j^- ) &= D_j \Big( [-\omega_0^*]_{-\omega^{M,min}_j}^{-\omega^{M,max}_j} \Big) \nonumber \\
&= [ D_j (-\omega_0^*)]_{p^{M,min}_j}^{p^{M,max}_j} \nonumber \\
&= p^{M,*}_j, \label{Dpeq}
\end{align}
by the above equivalences and equations~\eqref{eqbr5} and~\eqref{pf5eq}. Therefore, the optimality condition~\eqref{kkt1} holds. Analogously,
\begin{align}
D_{dj} ( \nu - \mu^+ + \mu^- ) &= D_{dj} \Big( [\omega_0^*]_{\omega^{c,min}_j}^{\omega^{c,max}_j} \Big) \nonumber \\
&= [ D_{dj} (\omega_0^*)]_{d^{c,min}_j}^{d^{c,max}_j} \nonumber \\
&= d^{c,*}_j. \label{Ddeq}
\end{align}
\end{subequations}
We now note that the preimages of any point $x \in \mathbb{R}$ under $D_j$ and $D_{dj}$ can be determined from the definitions~\eqref{geninverses} as
\begin{equation} \label{invgeninverses}
\begin{aligned}
D_j^{-1}(x) &= \begin{cases}
C_j' (x), &x \in \mathbb{R} \setminus \Lambda_j, \\
[C_j'(\gamma-),C_j'(\gamma+)], &x = \gamma \in \Lambda_j,
\end{cases} \\
D_{dj}^{-1}(x) &= \begin{cases}
C_{dj}' (x), &x \in {\color{black}\mathbb{R} \setminus \Lambda_{dj}}, \\
[C_{dj}'(\gamma-),C_{dj}'(\gamma+)], &x = \gamma \in \Lambda_{dj},
\end{cases}
\end{aligned}
\end{equation}
It then follows from~\eqref{invgeninverses}, the fact that $C_j'$ and $C_{dj}'$ are strictly increasing, and Definition~\ref{subdifdef} that, for all $x \in \mathbb{R}$,
\begin{equation} \label{subdifequalities}
D_j^{-1} (x) = \partial C_j (x) \text{ and } D_{dj}^{-1} (x) = \partial C_{dj} (x).
\end{equation}
Then, taking the preimages of~\eqref{Dpeq} and~\eqref{Ddeq} under $D_j$ and $D_{dj}$ respectively yields $-\nu - \lambda^+_j + \lambda^-_j \in D_j^{-1} (p^{M,*}_j)$ and $\nu - \mu^+_j + \mu^-_j \in D_{dj}^{-1} (d^{c,*}_j)$, whence we deduce by~\eqref{subdifequalities} that~\eqref{sdkkt1} and~\eqref{sdkkt2} are both satisfied by $\bar{p}^M_j = p^{M,*}_j$ and $\bar{d}^c_j =~d^{c,*}_j$.

Therefore, all of the subgradient KKT conditions~\eqref{sdkkt1},~\eqref{sdkkt2}, and~\eqref{kkt3}--\eqref{kkt8} are satisfied, and so $(p^{M,*}, d^{c,*}, d^{u,*})$ defines a global optimum for the OSLC problem~\eqref{Problem_To_Min}. \hfill \IEEEQED}

\balance

\begin{remark}
{\color{black}Theorem~\ref{Optimality_switch} shows that the optimality result given in Theorem~\ref{optthm} extends to situations in which the} cost functions are not required to be differentiable at all points. Such {\color{black}non-differentiable} cost functions {\color{black}can} describe {\color{black}important classes of} {\color{black}droop control schemes, such ones involving deadbands,}
{\color{black}as discussed in Section~\ref{Discussion}.}
\end{remark}


\begin{thebibliography}{31}

\bibitem{Schweppe}
F. C. Schweppe, R. D. Tabors, J. L. Kirtley, H. R. Outhred, F. H. Pickel, and A. J. {\color{black}Cox,} ``Homeostatic utility control,'' \emph{IEEE Trans. Power App.~Syst.}, vol.~PAS-99, no.~3, pp.~1151--1163, May~1980.

\bibitem{Recent Load Control}
J. A. Short, D. G. Infield, and L. L. Freris, ``Stabilization of grid frequency through dynamic demand control,'' \emph{IEEE Trans. Power Syst.}, vol.~22, no.~3, pp.~1284--1293, Aug.~2007.

\bibitem{Trudnowski}
D. Trudnowski, M. Donnelly, and E. Lightner, ``Power-system frequency and stability control using decentralized intelligent loads,'' IEEE PES T\&D Conf. Exhib., pp.~1453--1459, May~2006.

\bibitem{Primary2}
A. Molina-Garci\'{a}, F. Bouffard, and D. S. Kirschen, ``Decentralized demand-side contribution to primary frequency control,'' \emph{IEEE Trans. Power Syst.}, vol.~26, no.~1, pp.~411--419, May~2010.

\bibitem{Bergen_Vittal}
A. R. Bergen and V. Vittal, \emph{Power Systems Analysis}. Prentice Hall,~1999.

\bibitem{Pacific NW LAB}
D. J. Hammerstrom et al., ``Pacific Northwest GridWise testbed demonstration projects, part II: grid friendly appliance project,'' \emph{Pacific Northwest Nat. Lab.}, Tech. Rep. PNNL-17079, Oct.~2007.

\bibitem{LIPA}
B. J. Kirby, \emph{Spinning reserve from responsive loads}. United States Department of Energy,~2003.

{\color{black}\bibitem{P19}
C. Zhao, U. Topcu, N. Li, and S. Low, ``Design and stability of load-side primary frequency control in power systems,'' \emph{IEEE Trans. Autom. Control}, vol.~59, no.~5, pp.1177--1189, 2014.}

\bibitem{P34 - Low}
C. Zhao, and S. Low, ``Optimal decentralized primary frequency control in power networks,'' \emph{Proc. 53rd IEEE Conf. Decision Control}, pp 2467--2473, Dec.~2014.

\bibitem{LC2}
J. W. Simpson-Porco, F. D\"{o}rfler, and F. Bullo, ``Synchronization and power sharing for droop-controlled inverters in islanded microgrids,'' \emph{Automatica}, vol.~49, no.~9, pp.~2603--2611, Sep.~2013.

\bibitem{LC3}
Q. Shafiee, J. M. Guerrero, and J. C. Vasquez, ``Distributed secondary control for islanded microgrids -- a novel approach,'' \emph{IEEE Trans. Power Electron.}, vol.~29, {\color{black}no.~2}, pp.~1018--1031, Apr.~2013.

\bibitem{LC4}
M. Andreasson, H. Sandberg, D. V. Dimarogonas, and K. H. Johansson, ``Distributed integral action: stability analysis and frequency control of power systems,'' \emph{Proc. 51st IEEE Conf. Decision Control}, pp.~2077--2083, Dec.~2012.

\bibitem{trip}
S. Trip, M. B\"{u}rger, and C. De Persis. ``An internal model approach to frequency regulation in inverter-based microgrids with time-varying voltages,'' \emph{Proc. 53rd IEEE Conf. Decision Control}, pp 223--228, Dec.~2014.

{\color{black}\bibitem{P52}
X. Zhang and A. Papachristodoulou, ``A real-time control framework for smart power networks with star topology,'' \emph{Proc. American Control Conf.}, 2013.}

{\color{black}\bibitem{P39}
H. Miyagi and A. R. Bergen, ``Stability studies of multimachine power systems with the effects of automatic voltage regulators,'' \emph{IEEE Trans. Autom. Control}, vol.~31, no.~3, pp.~210--215, 1986.}


%
%
%
%

{\color{black}\bibitem{paper2}
E. Devane, A. Kasis and I. Lestas, ``Primary frequency regulation with load-side participation Part II: beyond passivity approaches,'' \emph{}{\color{black}
(under review for publication)
Feb. 2016.}}

{\color{black}\bibitem{Khalil}
H. K. Khalil and J. W. Grizzle. \emph{Nonlinear systems, vol.~3}, Prentice Hall New Jersey, 1996.}


\bibitem{sastry}
S. Sastry, \emph{Nonlinear systems: analysis, stability and control}. Springer-Verlag New York,~1999.

{\color{black}\bibitem{P50}
A. Molina-Garcia, F Bouffard, and D. S. Kirschen, ``Decentralized demand-side contribution to primary frequency control,'' \emph{IEEE Trans. Power Syst.}, vol.~26, no.~1, pp.~411--419, 2011.}

\bibitem{19-24}
G. Rogers, \emph{Power System Oscillations}. Kluwer Academic Publishers,~2000.

\bibitem{19-25}
K. W. Cheung, J. Chow, and G. Rogers, \emph{Power System Toolbox, v.~3.0}. Rensselaer Polytechnic Institute and Cherry Tree Scientific Software,~2009.

\bibitem{Boyd}
S. Boyd and L. Vandenberghe, \emph{Convex optimization}. Cambridge University Press,~2004.

{\color{black}\bibitem{Cauchy_Schwarz}
M. J. Steele. \emph{The Cauchy-Schwarz master class: an introduction to the art of mathematical inequalities}, Cambridge University Press, 2004.}

{\color{black}\bibitem{P28}
P. Moylan and D. Hill, ``Stability criteria for large-scale systems,'' \emph{IEEE Trans. Autom. Control}, vol.~23, no.~2, pp.~143--149, 1978.}

{\color{black}\bibitem{Feedback_Control_Doyle}
J. C. Doyle, B. A. Francis, and A. Tannenbaum. \emph{Feedback control theory, vol.~1}, Macmillan Publishing Company New York, 1992.}

{\color{black}\bibitem{Rockafellar}
R. T. Rockafellar. \emph{Convex analysis}. Princeton University Press, 1970.}

\end{thebibliography}
\end{document}